\newif\ifarxiv
\arxivtrue  

\ifarxiv
\documentclass[12pt]{article}
\usepackage[T1]{fontenc}   
\usepackage[utf8]{inputenc}
\usepackage{textcomp}      
\usepackage{lmodern}       

\usepackage{fullpage}

\usepackage{graphicx}
\usepackage{pifont}
\usepackage{pgfplots}

\usepackage{amsmath,amssymb,amsthm,enumitem,mathtools}
\renewcommand{\qedsymbol}{$\blacksquare$}

\usepackage{multirow}
\usepackage{algorithm}
\usepackage[noend]{algpseudocode} 

\usepackage{csvsimple}	
\usepackage{booktabs}   
\usepackage{adjustbox}  

\usepackage[colorlinks=true, citecolor=cyan, linkcolor=cyan, urlcolor=cyan]{hyperref}

\usepackage[xindy, acronyms]{glossaries}   
\makeglossaries

\newcommand{\snote}[1]{}
\renewcommand{\snote}[1]{\textcolor{magenta}{\textbf{\small [Stefan: #1]}}}

\newcommand{\bnote}[1]{}
\renewcommand{\bnote}[1]{\textcolor{blue}{\textbf{\small [ Bartolomeo #1]}}}

\newcommand{\eg}{{\it e.g.}}
\newcommand{\ie}{{\it i.e.}}

\newcommand{\ones}{\mathbf 1}
\newcommand{\reals}{{\mbox{\bf R}}}
\newcommand{\integers}{{\mbox{\bf Z}}}
\newcommand{\naturals}{{\mbox{\bf N}}}
\newcommand{\symm}{{\mbox{\bf S}}}  
\newcommand{\prob}{{\mathbf P}}
\newcommand{\distrib}{{\mathcal P}}
\newcommand{\identity}{I}
\newcommand{\tpose}{T}
\newcommand{\nullspace}{{\mathcal N}}
\newcommand{\range}{{\mathcal R}}
\newcommand{\Rank}{\mathop{\bf Rank}}
\newcommand{\Tr}{\mathop{\bf Tr}}
\newcommand{\diag}{\mathop{\bf diag}}
\newcommand{\lambdamax}{{\lambda_{\rm max}}}
\newcommand{\lambdamin}{\lambda_{\rm min}}
\newcommand{\Expect}{\mathop{\bf E{}}}
\newcommand{\Prob}{\mathop{\bf Prob}}
\newcommand{\Co}{{\mathop {\bf Co}}} 
\newcommand{\conv}{\operatorname{\bf conv}}
\newcommand{\dist}{\mathop{\bf dist{}}}
\newcommand{\argmin}{\mathop{\rm argmin}}

\newcommand{\epi}{\mathop{\bf epi}} 
\newcommand{\Vol}{\mathop{\bf vol}}
\newcommand{\dom}{\mathop{\bf dom}} 
\newcommand{\intr}{\mathop{\bf int}}

\newcommand{\sign}{\mathop{\bf sign}}
\newcommand{\round}{\mathop{\bf round}}
\newcommand{\define}{\coloneqq}

\newcommand{\children}[1]{c(#1)}

\newcommand{\param}{\theta}
\newcommand{\paramspace}{\Theta}
\newcommand{\paramdist}{\vartheta}
\newcommand{\vars}{x}
\newcommand{\lowervars}{y}
\newcommand{\covs}{x}
\newcommand{\predcost}{\hat{c}}
\newcommand{\predvars}{\hat{x}}
\newcommand{\predlowervars}{\hat{y}}
\newcommand{\truelayercost}{c}
\newcommand{\layer}{l}
\newcommand{\maxlayer}{L}
\newcommand{\varspace}{X}
\newcommand{\lowervarspace}{Y}
\newcommand{\bigvarspace}{W}
\newcommand{\hierarchy}{H}
\newcommand{\gspo}{\text{GSPO}+}
\newcommand{\lossname}{\text{D1}}
\newcommand{\losscut}{\text{D}}
\newcommand{\cost}{f}
\newcommand{\ubcost}{g}
\newcommand{\nn}{\hat{c}}
\newcommand{\nnparam}{w}
\newcommand{\confparam}{\gamma}
\newcommand{\dataset}{\mathcal{D}}
\newcommand{\calibration}{\mathcal{C}}
\newcommand{\eval}{\mathcal{E}}
\newcommand{\hypothesisclass}{\mathcal{H}}
\newcommand{\SUB}{\text{SUB}}
\newcommand{\rademacher}{\mathcal{R}}
\newcommand{\costbd}{C}
\newcommand{\avlossname}{\text{AV}}
\newcommand{\ext}{\textbf{\text{ext}}}
\newcommand{\ESUB}{\text{ESUB}}
\newcommand{\conff}{h}
\newcommand{\confnn}{\psi}
\newcommand{\ub}{u}
\newcommand{\lb}{l}
\newcommand{\test}{\mathcal{T}}

\newcommand{\cmark}{\ding{51}}%
\newcommand{\xmark}{\ding{55}}%

\newcommand{\dimension}[1]{\text{d}(#1)}

\newcommand{\T}{\text{TRUE}}
\newcommand{\F}{\text{FALSE}}

\newtheorem{theorem}{Theorem}[section]  
\newtheorem{lemma}{Lemma}[section]  

\newtheorem{assumption}{Assumption}[section]

\theoremstyle{definition}

\theoremstyle{remark}

\newcommand{\E}{\mathbf{E}}

\newacronym{LO}{LO}{linear optimization}
\newacronym{QO}{QO}{quadratic optimization}
\newacronym{MIQO}{MIQO}{mixed-integer quadratic optimization}
\newacronym{MIO}{MIO}{mixed-integer optimization}
\newacronym{MILO}{MILO}{mixed-integer linear optimization}
\newacronym{MINLO}{MINLO}{mixed-integer nonlinear optimization}

\newcommand{\transpose}{T}

\newcommand{\ld}{large decision problem}
\newcommand{\oa}{fast optimization algorithm}
\newcommand{\Ld}{Large decision problem}
\newcommand{\Oa}{Fast optimization algorithm}

\newcommand{\probct}{c(\theta)}
\newcommand{\probbt}{b(\theta)}
\newcommand{\probc}{c}
\newcommand{\probA}{A}
\newcommand{\probb}{b}
\newcommand{\pvar}{\Theta}
\newcommand{\algo}{\mathcal{A}}
\newcommand{\graveroracle}{\mathcal{G}}
\newcommand{\seedspace}{\mathcal{S}}
\newcommand{\risk}{R}
\newcommand{\erisk}{\hat{R}}
\newcommand{\seed}{S}
\newcommand{\complexity}{\mathcal{K}}
\newcommand{\wcl}{\omega}

\newcommand{\controltime}{t}
\newcommand{\controlmaxtime}{T}
\newcommand{\power}{P}
\newcommand{\maxpower}{\power^{\text{max}}}
\newcommand{\samplingtime}{\upsilon}
\newcommand{\powerload}{\power^{\text{load}}}
\newcommand{\energy}{E}
\newcommand{\minenergy}{\energy^{\text{min}}}
\newcommand{\maxenergy}{\energy^{\text{max}}}
\newcommand{\initenergy}{\energy^{\text{init}}}
\newcommand{\fuelstate}{\zeta}
\newcommand{\alphacontrol}{\alpha}
\newcommand{\betacontrol}{\omega}
\newcommand{\gammacontrol}{\delta}
\newcommand{\dcontrol}{\xi}
\newcommand{\hcontrol}{h}
\newcommand{\wcontrol}{\psi}
\newcommand{\fuelstateinit}{\fuelstate^{\text{init}}}
\newcommand{\scontrol}{k}
\newcommand{\scontrolinit}{s^{\text{init}}}
\newcommand{\numswitchingsinit}{\numswitchings_{\text{init}}}
\newcommand{\dcontrolpast}{\dcontrol^{\text{past}}}
\newcommand{\nsw}{N^{\text{sw}}}
\newcommand{\Gcontrol}{G}
\newcommand{\controlf}{f}
\newcommand{\dpast}{\dcontrol^{\text{past}}}
\newcommand{\weightset}{\mathcal{W}}
\newcommand{\indicator}{\mathbf{1}}
\newcommand{\ASL}{\text{ASL}}
\newcommand{\MASL}{\text{MASL}}

\newcommand{\numconstraints}{p}
\newcommand{\numvars}{n}

\newcommand{\numhouses}{H}
\newcommand{\numroutes}{R}
\newcommand{\numrouteconstraints}{p}

\newcommand{\no}{n_1}
\newcommand{\nt}{n_2}
\newcommand{\ncon}{m}
\newcommand{\paramdim}{p}
\newcommand{\costx}{c}
\newcommand{\costy}{d}
\newcommand{\demand}{\delta}
\newcommand{\HO}{\text{D}}
\newcommand{\GRB}{GRB}
\newcommand{\GRBo}{GRB-H}
\newcommand{\GRBt}{GRB-3SOL}
\newcommand{\SCIP}{SCIP}
\newcommand{\SCIPo}{SCIP-H}
\newcommand{\SCIPt}{SCIP-3SOL}

\newcommand{\TIMELIMIT}{100}

\newcommand{\change}[1]{#1}
\newcommand{\NINSTANCE}{10000}
\newcommand{\NTRAIN}{9700}
\newcommand{\NEVAL}{100}
\newcommand{\NCAL}{100}
\newcommand{\NTEST}{100}
\newcommand{\NTHETA}{100}
\newcommand{\NLAYER}{3}
\newcommand{\NNEURON}{2000}

\newcommand{\NKNAPSACK}{100}

\newcommand{\NFACLOCCOMPLICATING}{25}
\newcommand{\FACLOCGAMMA}{100}
\newcommand{\FACLOCN}{75}

\newcommand{\ROUTINGNNODE}{150}
\newcommand{\ROUTINGNVEHICLE}{20}
\newcommand{\ROUTINGCAPACITY}{10}
\newcommand{\ROUTINGFUEL}{10}

\newcommand{\STOVSPNTASK}{50}
\newcommand{\STOVSPM}{50}

\newcommand{\TOLERANCE}{0.0001}

\newcommand{\FY}{\text{FY}}
\newcommand{\fenchel}{\Omega}

\newcommand{\GEN}{\text{GEN}}
\newcommand{\LIN}{\text{LIN}}
\newcommand{\smallloss}{\ell}

\newcommand{\CPALPHA}{0.1}

\newcommand{\faclocz}{\eta}
\newcommand{\routingz}{s}
\newcommand{\routingw}{\omega}
\newcommand{\stovspslack}{\eta}

\newcommand{\cstmcite}[1]{\cite{#1}}
\newcommand{\cstmciter}[1]{\cite{#1}}

\else
\documentclass[ijoo,dblanonrev]{informs4}  
\usepackage{eqndefns-left} 
\RequirePackage{tgtermes}
\RequirePackage{newtxtext}
\RequirePackage{newtxmath}
\RequirePackage{bm}
\RequirePackage{endnotes}

\usepackage{booktabs}   
\usepackage{csvsimple}  
\usepackage{pifont}     

\OneAndAHalfSpacedXII 

\usepackage{algorithm}
\usepackage{algpseudocode}
\usepackage{tikz}

\EquationsNumberedThrough    

\TheoremsNumberedThrough     
\ECRepeatTheorems  %

\MANUSCRIPTNO{IJOO-0001-2024.00}

\newcommand{\snote}[1]{}
\renewcommand{\snote}[1]{\textcolor{magenta}{\textbf{\small [Stefan: #1]}}}

\newcommand{\bnote}[1]{}
\renewcommand{\bnote}[1]{\textcolor{blue}{\textbf{\small [ Bartolomeo #1]}}}

\newcommand{\cstmcite}[1]{\citep{#1}}
\newcommand{\cstmciter}[1]{\citet{#1}}

\fi

\title{Learning-Based Hierarchical Approach for \\Fast Mixed-Integer Optimization}
\author{Stefan Clarke and Bartolomeo Stellato}
\date{\small Department of Operations Research and Financial Engineering\\[.5em] Princeton University\\[.5em]\today}

\begin{document}

\ifarxiv
\maketitle

\begin{abstract}
    We propose a hierarchical architecture for efficiently computing high-quality solutions to structured mixed-integer programs (MIPs). 
    To reduce computational effort, our approach decouples the original problem into a higher level problem and a lower level problem, both of smaller size.
    We solve both problems sequentially, where decisions of the higher level problem become parameters of the constraints of the lower level problem.
    We formulate this learning task as a convex optimization problem using decision-focused learning techniques and solve it by differentiating through the higher and the lower level problems in our architecture.
    To ensure robustness, we derive out-of-sample performance guarantees using conformal prediction.
    Numerical experiments in facility location, knapsack problems, and vehicle routing problems demonstrate that our approach significantly reduces computation time while maintaining feasibility and high solution quality compared to state-of-the-art solvers.
\end{abstract}

\fi

\section{Introduction}
Many decision problems in engineering, computer science, and operations research, involve repeatedly solving similar optimization problems with varying data, \ie, \emph{reoptimizing}.
For example, in power-grid operations, the structure of the grid remains fixed while certain parameters, such as the power demand and the wind forecasts vary~\cstmcite{pascal, Zamzam2019LearningOS, Bienstock2020}.
Routing problems can also exhibit a similar structure where the road network is fixed while the demand for each location may vary~\cstmcite{VRPBook, Bertsimas2019OnlineVR, Jaillet2006OnlineRP, Bose1999OnlineRI}. 
In these applications, instead of resolving each problem from scratch, we may be able to reuse the information from previous instances to accelerate subsequent solutions.

Due to recent hardware and software progress state-of-the-art mixed-integer programming (MIP) solvers are able to solve MIPs at remarkable speeds, despite the fact that they are NP-hard~\cstmcite{Bixby2007, Avella2023ACS}.
However, their execution time may not be fast enough in many real-time applications. 
For example, problems arising in the optimization of power systems often need to be solved in minutes~\cstmcite{Namor2018ControlOB}, and control problems in robotics may need to be solved in milliseconds~\cstmcite{Liniger2015OptimizationbasedAR, Raibert2008BigDogTR}.

Many industrial settings provide access to significant amounts of data about problem parameters and historical realizations.
In robotics, for example, engineers often have recorded previous states visited by the robot or typical environmental parameters.
This data can be used to train models that accelerate optimization.
The simplest approach is to train a model that directly predicts the optimization solution, \ie, the end-to-end method in~\cstmcite{Bengio2018MachineLF}.
However, such predictions come with no guarantees on feasibility or optimality.
More recent works exploit data differently: using learning to guide the solution process for mixed-integer programs, which allows real-time computation of high-quality solutions~\cstmcite{Zhang2022ASF, Bengio2018MachineLF, Alvarez2017AML}.
We contribute to this line of work by developing learned optimizers for structured mixed-integer optimization problems.

A key aspect of many applications is that fixing certain variables decomposes a challenging optimization problem into smaller, more tractable subproblems.
For example, in multi-agent routing, once we fix the target locations for each agent, each agent can independently plan its trajectory, yielding smaller problems that we can solve more quickly~\cstmcite{Zhang2020MultiVehicleRP, Silva2019ARL, Braysy2005, Desrochers1992}.
In facility location problems, once we fix the facility allocations, the problem becomes continuous~\cstmcite{Liu2009, Avella2023ACS}.
After fixing these top-layer variables, the problem simplifies considerably: instead of jointly locating facilities and assigning customers, we only need to determine the continuous flow variables between facilities and customers.
Such key decisions, often referred to as \emph{backdoors}~\cstmcite{Cai2024LearningBF, Fischetti2011BackdoorB, bistra}, once found, greatly simplify the solution of the entire problem.  
In these cases, we can exploit the problem structure and solve it hierarchically: first optimize over the most important variables, the higher level problem, and then solve for the remaining variables, the lower level problem.
Since the higher level problem depends only on a subset of the decision variables, it is smaller than the original MIP and faster to solve.
This approach typically does not yield the optimal solution to the original problem, and a poor choice of higher level variables may even render the lower level problem infeasible.
The design of such a hierarchical decomposition therefore requires careful attention.

In this paper, we focus solving hierarchically structured MIPs by learning small, tractable formulations consisting of an upper level and a lower level problem, which can be efficiently solved to obtain heuristic solutions.
Our contributions are as follows.
\begin{itemize}
\item We propose a method to learn to approximately solve parametric mixed-integer programs.
We introduce a differentiable two-layer architecture that computes a candidate solution by first solving a high level problem, then using its solution to parametrize a lower level problem whose solution yields the final candidate point.
The key component of this architecture is a neural network that predicts the coefficients of the high level objective function to optimize the quality of the candidate solution.

\item We formulate the training problem as a convex optimization problem that a convex loss, inspired by the Smart Predict and Optimize~\cstmcite{spo} and the Fenchel-Young~\cstmcite{Dalle2022LearningWC} losses.
While the training problem is convex, evaluating the loss repeatedly over many iterations is challenging since each evaluation requires solving a potentially large mixed-integer program.
We therefore develop several surrogate convex losses that approximate our generic loss, and provide meaningful gradients, and have the same minimizers.

\item We develop a method to obtain probabilistic bounds on the suboptimality of predicted solutions using conformal prediction techniques~\cstmcite{conformal,conformal_vovk}. This method trains a neural network on an evaluation dataset to predict the optimal value, then uses a separate calibration dataset to construct validity guarantees via conformal prediction.

\item We demonstrate the effectiveness of our methods with a series of computational experiments on randomly generated problem instances from facility location, knapsack, and vehicle routing problems.
We show that our approach is able to find feasible solutions significantly faster than state-of-the-art MIP solvers while maintaining high solution quality.


\end{itemize}

\subsection{Related literature}

\paragraph{Learning for optimization.} Machine learning has been used to speed up optimization algorithms in many settings~\cstmcite{Amos2022TutorialOA, learntooptimize}, including in discrete optimization~\cstmcite{Bengio2018MachineLF}, power grid optimization~\cstmcite{pascalenergy}, control problems~\cstmcite{Mitrai2024AcceleratingPC}, and finance~\cstmcite{Rychener2023EndtoEndLF}.
Some attempts to learn to optimize combinatorial problems focus on directly modifying the solution algorithm, \eg, branch-and-bound~\cstmcite{Scavuzzo2024MachineLA}, while some focus on using machine learning to directly guess the values of integer variables~\cstmcite{Bertsimas2018TheVO}.
Lack of differentiability of the solution with respect to parameters makes the problem of learning in MIPs difficult.
Also learning algorithms often require the repeated solution of costly subproblems.
Unlike other works, we will attempt to formulate the problem of learning a fast hierarchical integer optimization solver into a convex problem and in such a way that solving subproblems is not too costly.

\paragraph{Decision-focused learning and inverse optimization.}
Decision-focused learning (DFL)~\cstmcite{Mandi2023DecisionFocusedLF} is a framework for training predictors that optimize the quality of downstream optimization tasks rather than minimizing prediction error.
In this setting, the predictor outputs become parameters of an optimization problem, and the goal is to learn predictions that lead to high-quality solutions~\cstmcite{Wilder2018MeldingTD}.
Decision focused learning has applications in several areas, including reinforcement learning~\cstmcite{Sharma2023DecisionFocusedMR}, resource allocation problems~\cstmcite{Verma2022CaseSA}, and finance~\cstmcite{Lee2024AnatomyOM}.
When learning to predict the objective vector of a MIP, the supervised DFL task can be reformulated into a convex optimization problem~\cstmcite{spo}.
In this work, rather than using DFL for prediction and regression tasks, as is currently common in the literature, we use it to learn faster sequential optimizers for large optimization problems.

\paragraph{Differentiable learning of mixed integer programs.}
Some recent works have attempted to exploit differentiability in specific ways to machine-learn integer programming solvers~\cstmcite{diffcp, ShiviGupta}. These works have mostly been focused on generation of cutting planes because there are many ways to generate cutting planes in a differentiable way. Other works have formulated MIPs in an approximately differentiable way so that they can be used as layers in machine learning architectures~\cstmcite{Ferber2019MIPaaLMI}.
Some have used rounding and probabilistic techniques to approximate MIPs with differentiable problems~\cstmcite{geng2025differentiable}. In each of these cases the learning problem is rarely convex. In this work we create a convex reformulation of a machine-learning task for integer optimization which involves the learning of the objective vector in a convex way, and which can be used to train small and simple models to perform in complex tasks. 
    Previous work has been done in using decision--focused--learning--like techniques to attempt to solve intractible optimization problems quickly. 
    Dalle et. al.~\cstmcite{Dalle2022LearningWC} use Fenchel-Young loss functions~\cstmcite{Blondel2019LearningWF} to approximate intractible optimization problems with tractable ones so that they can solve the tractable problem rather than solving the larger one, and extract a feasible solution for the harder problem from one for the easier one.
    They also provide methods for learning to solve these optimization problems in an unsupervised way. 
    In this paper we take a similar approach. 
    The loss function surrogates we consider in this paper can be interpreted as versions of the Support Vector Machine (SVM) loss~\cstmcite{Dalle2022LearningWC}.
    We use various modifications of the Fenchel-Young loss function to learn to solve integer optimization problems in a hierarchical way, and also derive suboptimality bounds for our approximations online.   


\subsection{A motivational example: a production problem} \label{sec:product}
The motivation and setup for our work can be explained through an inventory optimization example of a product supplier company.
In a given month the company receives orders for~$p \in \integers_+$ products in total. We assume that the amount of products which must be produced $p$ vary through some parameter space $\paramspace$. 
We assume that the company has some historical data of past realizations of~$p$. 
The company owns $f$ factories, 
and needs to supply each with the resources necessary for production~$s \in \integers_+$.
Given $r^1_i, \dots, r^s_i \in \reals$ resources from the company and an order of $\bar{d}_i$ products, the production plan for each factory $i \in \{1, \dots f\}$ solves the following optimization problem,
\begin{equation}
\begin{array}{ll} \label{eq:fac1} \ifarxiv \tag{$\mathcal{P}_{\text{F}}^i$} \fi
    \textstyle
    \text{minimize} 
     & {c_i}^\transpose \faclocz_i  + p_i (\bar{d}_i - d_i)_+   \\
         \text{subject to}   & A_i \faclocz_i + a_i d_i \le h(r_i),\\
         & d_i \ge 0, \quad d_i \in \integers, \quad \faclocz_i \in \integers^{n_i},
\end{array}
\end{equation}
where the decision variable $\faclocz_i \in \integers^{n_i}$ represents all decisions made by the factory (\eg, which machines are operated, how many workers are hired), and $d_i \in \integers$ the number of products produced by factory $i$ in total. 
The parameters are matrix~$A_i \in \reals^{m \times n_i}$, a vector valued function~$h:\reals \mapsto \reals^m$ which depends on the amount of resources $r_i=(r^1_i, \dots r^s_i) \in \reals^s$ assigned to factory $i$, and the vector $c_i \in \reals^{n_i}$ representing the cost of production. 
Then the variables for the production problem for factory~$i$ are given by $y_i = (\faclocz_i, d_i)$. 
The factory problems can be combined into a single optimization problem with variable $y=(y_1, \dots y_q)$ which is separable into each of the factories.

\begin{figure}
    \centering \includegraphics[width=0.9\linewidth, trim={0cm 18cm 0 2cm}]{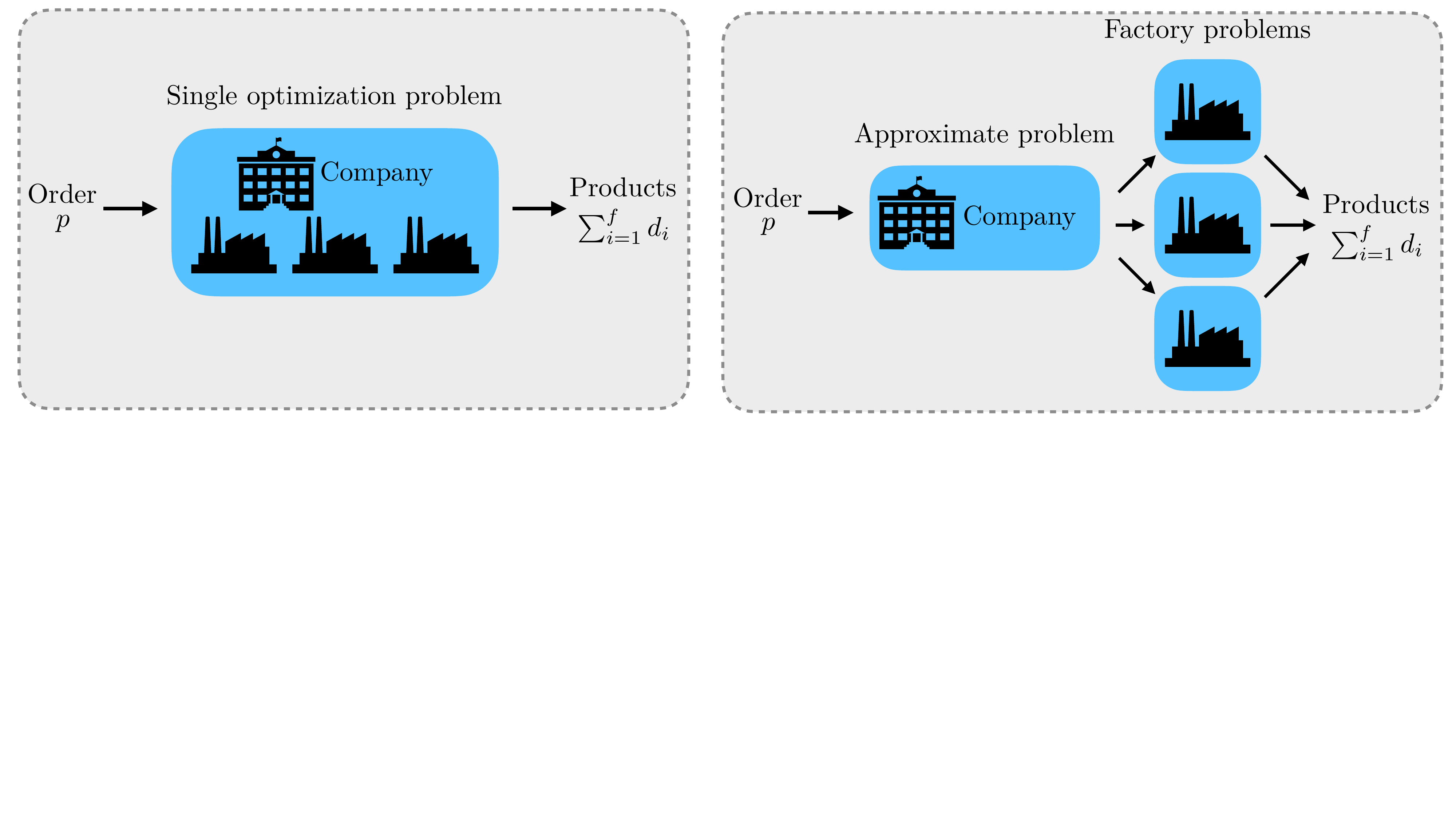}
    \caption{Top: Single, large optimization problem. Exact but may be intractible. Bottom: Hierarchical integer programming problems to approximate the top problem. Inexact but fast, with learnable parameters.
    }
    \label{fig:product}
\end{figure}

The company wants to avoid missing orders while also minimizing total costs from the factories, the purchasing and delivery of resources, and the delivery of the produced products to the warehouse. Overall the company is interested in solving the following problem,
\begin{equation} \label{eq: bigproduct} \ifarxiv \tag{$\mathcal{P}_{\text{F}}$} \fi
    \begin{array}{ll}
        \text{minimize}  & \textstyle g (p - \sum_{i=1}^f d_i)_+ + \sum_{i=1}^f c_i^\transpose \faclocz_i + {e_i} d_i + h^\transpose r_i \\
         \text{subject to} & \textstyle A r \le b, \\ & A_i \faclocz_i + a_i d_i \le h(r_i),\\
         & \textstyle r_i \in \integers^{s} \quad i=1, \dots f, \quad r = (r_1, \dots r_f),
    \end{array}
\end{equation}
where $e_i \in \reals$ is the delivery cost for products from factory $i$, $g \in \reals$ is the penalty for missed orders, $h \in \reals^s$ is the cost of each resource, and the constraint $A r \le b$ represents a polyhedral constraint on the number of resources it is possible to purchase. 
The corresponding dimensions are~$A \in \reals^{q \times (fs)}$ and $b \in \reals^q$ for some $q \in \integers_+$.
If the number of factories~$f$ is very large, solving \eqref{eq: bigproduct} exactly might be computationally expensive and therefore impractical for the company.
In this case, the company may choose to solve a smaller problem to decide on the amount to order from each factory, send these numbers to the factories, and then let each factory independently solve its own problem of the form \eqref{eq:fac1} and send the products back. 
This would mean sacrificing guaranteed optimality of the integer program in exchange for a shorter solve time. In this case, the company would purchase resources and make orders based on the following optimization problem,
\begin{equation} \label{eq:smallproduct}
    \textstyle
    \begin{array}{ll}
        \text{minimize} \quad & \textstyle \hat{g} (p - \sum_{i=1}^f \bar{d}_i)_+ + \sum_{i=1}^f  {\hat{e}_i} \bar{d}_i + \hat{h}_i^\transpose r_i  \\
         \text{subject to} \quad & A r \le b, \\
         & \textstyle \bar{d}_i \in \integers,\quad  r_i \in \integers^{s}, \\ & \textstyle \bar{d} = (\bar{d}_1, \dots \bar{d}_f), \quad r = (r_1, \dots r_f),
    \end{array} \ifarxiv \tag{$\hat{\mathcal{P}}_{\text{F}}$} \fi
\end{equation}
where $\hat{g} \in \reals, \hat{e}_i \in \reals,$ and $ \hat{h}_i \in \reals^s$ have been chosen so that the problem \eqref{eq:smallproduct} accurately approximates \eqref{eq: bigproduct}. This creates a hierarchical system of integer optimization problems, as shown in Figure~\ref{fig:product}. The optimization variable for the upper problem in the hierarchy is given by $\vars = (\bar{d}, r)$. 
It is not obvious exactly how to choose the parameters $\hat{g}, \hat{e}$ and $\hat{h}$ so that the hierarchical model accurately mirrors the full model.
In this paper, we focus on learning these terms in the objective so that the problems~\eqref{eq: bigproduct} and~\eqref{eq:smallproduct} have optimal solutions with similar cost. 
While the resulting solutions of~\eqref{eq:smallproduct} given by this process are guaranteed to be feasible, there is no guarantee that they will be optimal in~\eqref{eq: bigproduct}. A natural question for the company owners to ask is, \emph{how close to optimal is the solution given by~\eqref{eq:smallproduct}?} Had the company chosen to solve~\eqref{eq: bigproduct} using a branch-and-bound solver and terminated the solve early, the company would have this information in the form of a lower-bound for the objective. However a bound for the objective of~\eqref{eq: bigproduct} does not naturally lead to a bound for the objective of~\eqref{eq:smallproduct}. In Section~\ref{sec:lowerbounds} we will formulate methods for deriving rigorous bounds on the optimality of the solutions to problem~\eqref{eq:smallproduct} in problem~\eqref{eq: bigproduct}.

\section{Problem setup} \label{sec:setup}
Let $\paramspace \subseteq \reals^\paramdim$ be a set and $\param \in \paramspace$ be a random vector distributed according to distribution $\paramdist$ supported on $\paramspace$.
Let the decision of the \emph{upper problem} be~$x \in \varspace(\param) \subseteq \reals^{\no}$, and the decisions of the \emph{lower problem} by $y \in \lowervarspace(\vars, \param) \subseteq \reals^{\nt}$. 
Let $\costx(\param) \in \reals^{\no}$ and $\costy(\param) \in \reals^{\nt}$ be the objective coefficients of the upper and lower problems respectively.
Our goal is to solve the following optimization problem over variables $x$ and $y$,
\begin{equation} \label{eq:master} \ifarxiv \tag{$\mathcal{P}$} \fi
    z(\param) = \begin{array}[t]{ll}
         \text{minimize}  & {\costx(\param)}^\transpose \vars + {\costy(\param)}^\transpose \lowervars \\
         \text{subject to} & \vars \in \varspace(\param), \quad \lowervars \in \lowervarspace(\vars, \param).
    \end{array}
\end{equation}
In many examples, such as in~\eqref{eq: bigproduct}, once the variable in the upper problem problem $\vars$ is fixed the problem~\eqref{eq:master} might become separable, meaning that finding different components of $\lowervars$ can be done in parallel.
In such examples the problem can be rewritten as,
\begin{equation} \label{eq:decomp}
    \begin{array}{ll}
         \text{minimize}  & c(\param)^\transpose \vars + \sum_{i=1}^{k}d_i(\param)^\transpose \lowervars_i \\
         \text{subject to} & \vars \in \varspace(\param), \quad \lowervars_i \in \lowervarspace_i(\vars, \param)\quad i=1,\dots,k.
    \end{array}
\end{equation}
We make the following key assumption of the feasible region.
\begin{assumption}[Recursive feasibility]\label{ass:recursive}
    The decision set $\lowervarspace(\vars, \param)$ is nonempty for every $\vars \in \varspace(\param)$ and $\param \in \paramspace$.
\end{assumption}
This ensures that there is always a feasible decision in the lower \change{problem}---even with an arbitrarily high cost---no matter what decision we make at the top layer.
Many discrete optimization problems have this structure, such as the problems we consider in Section~\ref{sec:computation}. 

\subsection{Learning an optimizer}
We aim to solve the problem~\eqref{eq:master} by first deciding $\vars$ and then $\lowervars$, so that instead of solving a single large problem we can solve a sequence of smaller ones, with the intention of decreasing solve-time.
Define the true cost of making decision $\vars$ in the first layer as,
\begin{equation} \label{eq:truecost}
    \cost_\param(\vars) = \textstyle c(\param)^\transpose \vars + \begin{array}[t]{ll}
        \text{minimize}  &  \textstyle d(\param)^\transpose \lowervars \\
                         &\text{subject to} \quad \lowervars \in \lowervarspace(\vars, \param).
    \end{array}
\end{equation}
The function $\cost$ is hard to evaluate because it requires solving a hard optimization problem to compute the lower level decisions~$y$, so we approximate it with a linear function~$\predcost(\param)^Tx$. 
Therefore, we learn methods which make a decision for the upper-level problem by solving,
\begin{equation} \label{eq:costform1} \ifarxiv \tag{$\hat{\mathcal{P}}^{1}$} \fi
    \hat{\vars}(\param) = \begin{array}[t]{ll}
        \argmin & \predcost(\param)^\transpose \vars 
        \qquad  \text{subject to} \quad \vars \in \varspace(\param).
    \end{array}
\end{equation}
We refer to this problem as our \emph{policy} as the resulting decisions $\hat{\vars}(\param)$ are the core part of our problem.
Given $x$, the lower decisions can be computed as
\begin{equation} \label{eq:lowerprob}\ifarxiv \tag{$\mathcal{P}^{2}$} \fi
    \hat{\lowervars} = \begin{array}[t]{ll}
        \argmin  & \costy(\param)^\transpose \lowervars \qquad
         \text{subject to} \quad \lowervars \in \lowervarspace(\vars, \param).
    \end{array}
\end{equation}
\begin{figure}
    \centering
    \includegraphics[trim={0 25cm 0 0cm},clip,width=0.8\linewidth]{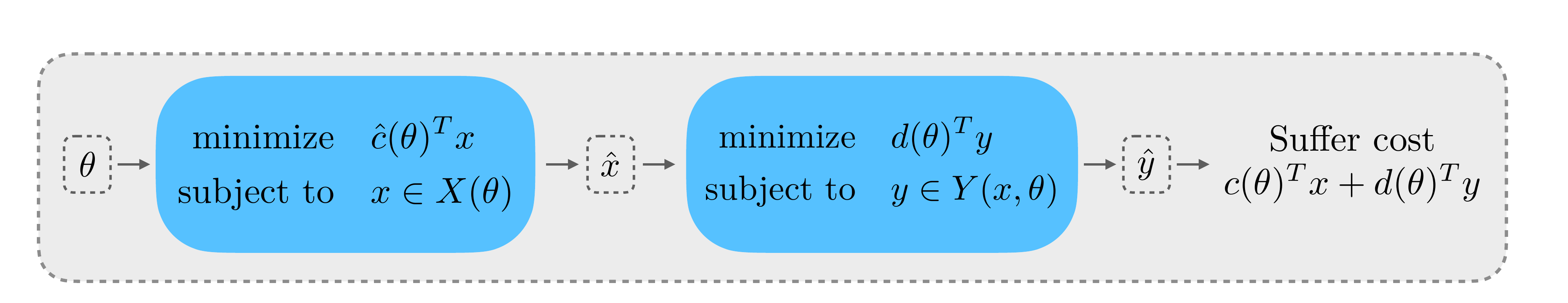}
    \caption{The layered optimization model setup used to compute feasible solutions and upper-bounds.}
    \label{fig:layeredsetup}
\end{figure}
Our goal is to learn a function $\predcost(\param)$ to parametrize problem~\eqref{eq:costform1} so that its solutions, and the corresponding lower level solutions to~\eqref{eq:lowerprob}, are feasible and achieve good objective values for the original problem~\eqref{eq:master}.
We model the predictor~$\predcost$ as a function with weights~$\nnparam$, \eg, a neural network.
We represent this learning task as minimizing the \emph{suboptimality loss}~\cstmcite{spo},
\begin{equation} \label{eq:subloss}
     \mathcal{L}_{\text{\SUB}}(\nnparam) = \underset{\param \sim \paramdist}{\E} \left[ \smallloss_{\text{SUB}}(\predcost_\nnparam(\param) | \param)\right],\quad\text{with}\quad  \smallloss_{\text{SUB}}(\costx | \param) =  \underset{\vars \in \varspace^\star(\costx | \param)}{\max} \cost_\param(\vars) 
     - \cost_\param(\vars^\star(\param)),
\end{equation}
where, for any $\param$, $\varspace^\star(\costx | \param)$ is the set of minimizers of $\costx^\transpose \vars$ on $\varspace(\param)$ and $(\vars^\star(\param),  \lowervars^\star(\param))$ solve~\eqref{eq:master} for a given $\param \in \paramspace$.
Informally, for a given $\param$ and $\predcost$, the suboptimality loss can interpreted as penalizing $\costx(\param)^\transpose \big(\predvars - \vars^{\star}(\param)\big) + \costy(\param)^\transpose \big(\predlowervars - \lowervars^{\star}(\param)\big)$.
However, for a given $\predcost$ there could be multiple optimal solutions $\predvars$. 
That's why we adopt the more precise definition~\eqref{eq:subloss}.
The loss function~\eqref{eq:subloss} is often nonconvex, discontinuous, and even piecewise constant. Therefore optimizing directly over parameters in $\predcost$ is very difficult. To remedy this issue, we propose convex surrogate loss functions for~\eqref{eq:subloss} and give methods to optimize them, as well as theoretical guarantees that the optima of the true suboptimality loss and the surrogates are close.

\section{Differentiable convex surrogate functions for the suboptimality loss} \label{sec:surrogate}
In this section, we consider the loss function for a single value of $\param$ and, for notational convenience, we drop the dependency on the instance parameter $\theta$ from $\vars^\star$, $\smallloss_{\SUB}$, and~$\predcost$, and the dependency on the predictor weights~$\nnparam$ from~$\predcost$.
Our goal is to find $\predcost$ which minimizes $\smallloss_{\SUB}(\predcost)$.
Let $\ubcost : \varspace \mapsto \reals$ be any function satisfying the following assumptions.
\begin{assumption} 
    \label{ass:cost}
        The function $\ubcost(\vars)$ attains its minimum at $\vars^\star$ on $\varspace$,
        is bounded above and below on $\varspace$,
        and is upper semicontinuous.
\end{assumption}
We consider \emph{general} loss functions which approximate~\eqref{eq:subloss} with the form,
\begin{equation} \label{eq:linloss} \ifarxiv \tag{$\LIN$} \fi
    \smallloss_{\GEN}(\predcost) = \underset{\vars \in \varspace}{\max} \Big\{ \ubcost(\vars) - \ubcost(\vars^\star) -  \predcost^\transpose (\vars - \vars^\star) \Big\}.
\end{equation}

The following Theorem shows that these loss functions are useful surrogates because their optima line up with the true optima of~\eqref{eq:subloss}, they are convex, and their subgradients are easy to compute. 
This is inspired by known results on the SPO+~\cstmcite{spo} and the SVM~\cstmcite{Dalle2022LearningWC} losses.
\begin{theorem} \label{thm:derivative}
    The following are true for any $\ubcost$ satisfying Assumption~\ref{ass:cost}. 
    \begin{enumerate}
        \item $\smallloss_{\GEN}(\predcost)$ is a convex function. \label{thm:derivative:2}
        \item A subgradient of $\smallloss_{\GEN}(\predcost)$ at $\predcost$ is given by $(\vars^\star - \vars)$ where $\vars$ solves the optimization problem in \eqref{eq:linloss}.
            More generally, if $\vars_{\epsilon}$ is $\epsilon$-optimal in \eqref{eq:linloss} then $(\vars^{\star}- \vars_\epsilon )$ is an $\epsilon$-subgradient of $\smallloss_{\GEN}(\predcost)$. \label{thm:derivative:3}
        \item If $\smallloss_{\GEN}(\predcost) = 0$ then also $\smallloss_{\SUB}(\predcost) = 0$ and so any $\vars$ which minimizes $\vars \mapsto \predcost^\transpose \vars$ on $\varspace$ is optimal for $\ubcost$.
            Further, if $\vars^\star$ is the unique minimizer of $\ubcost$ on $\varspace$ then $\vars^\star$ is also the unique minimizer of $\predcost^\transpose \vars$ on $\varspace$. \label{thm:derivative:3.5}
        \item If $\vars^\star$ is an extreme point of $\varspace$ then for any ${\epsilon > 0}$ there is a $\predcost \in \reals^n$ such that ${\smallloss_{\GEN}(\predcost) < \epsilon}$.\label{thm:derivative:5}
    \end{enumerate}
\end{theorem}

In order to prove Theorem~\ref{thm:derivative}, we make use of the following lemma. 
\begin{lemma}
    \label{lemma:sharpmin} Let $\varspace \subseteq \reals^n$ be a compact set, $\ubcost$ be a function satisfying Assumption~\ref{ass:cost}, and $\vars^\star$ be an extreme point of $\varspace$. 
    Then, for any $\epsilon > 0$ there exists a vector $\truelayercost \in \reals^n$ such that $\vars^\star$ minimizes $\vars \mapsto \truelayercost^\transpose \vars$ on $\varspace$
    and such that for any $\vars \in \varspace$, $\truelayercost^\transpose (\vars - \vars^\star) \ge \ubcost(\vars) - \ubcost(\vars^\star) - \epsilon$.
\end{lemma}
\ifarxiv
    \begin{proof}
    \else
        \proof{Proof of Lemma~\ref{lemma:sharpmin}.}
    \fi
    Because $\vars^\star$ is an extreme point of $\varspace$, we can choose $\tilde{\costx} \in \reals^n$ such that $\vars^\star$ is the unique minimizer of  $\vars \mapsto \tilde{\costx}^\transpose \vars$ on $\varspace$. For $\delta > 0$ define the halfspace $R_\delta = \{ \vars \mid \tilde{\costx}^\transpose (\vars - \vars^\star) \le \delta\}$ and $S_\delta = \varspace \backslash R_\delta$. Let $C$ be an upper bound for $\ubcost$ on $\varspace$ and $\lambda > C/\delta$. Then, for any $\vars \in S_\delta$ we have
    $\lambda \tilde{\costx}^\transpose (\vars - \vars^\star) \ge (C / \delta) \delta = C \ge \ubcost(\vars).$
    Since $\vars^\star$ is the unique minimizer of $\vars \mapsto \tilde{\costx}^\transpose \vars$ on $\varspace$, as $\delta \rightarrow 0$ the set $R_\delta$ shrinks to $\{\vars^\star\}$.
    Therefore by upper semicontinuity of $\ubcost$, we can choose $\delta$ to be such that for all $\vars \in R_{\delta},$ $\ubcost(\vars) - \ubcost(\vars^\star) \le \epsilon$. Then, for any $\vars \in S_\delta$,
    $\lambda \tilde{\costx}^\transpose (\vars - \vars^\star) \ge 0 \ge \ubcost(\vars) - \ubcost(\vars^\star) - \epsilon.$
    Therefore we can set $\costx = \lambda \tilde{\costx}$ for the chosen values of $\delta$ and $\lambda$ and we have for $\vars$ both in $R_\delta$ and $S_\delta$ that $\costx^\transpose (\vars - \vars^\star) \ge \ubcost(\vars) - \ubcost(\vars^\star) - \epsilon$.
\ifarxiv
\end{proof}
\else
\qedsymbol
\endproof
\fi

We are now ready to prove Theorem~\ref{thm:derivative}.
\ifarxiv
    \begin{proof}
    \else
        \proof{Proof of Theorem~\ref{thm:derivative}.}
    \fi
    For item~\ref{thm:derivative:2}, observe that $\smallloss_{\GEN}$ is a maximum of convex functions of $\predcost$ and, therefore, convex.
    For item~\ref{thm:derivative:3}, for any $\tilde{c} \in \reals^n$ we have,
    \begin{align*}
        \textstyle \smallloss_{\GEN}(\predcost) + (\vars^\star - \vars_\epsilon)^\transpose (\tilde{c} - \predcost) &= \underset{\vars \in \varspace}{\max}\{ \ubcost(\vars )
        - \predcost^\transpose (\vars - \vars_\epsilon )\} - \tilde{c}^\transpose \vars_\epsilon + \tilde{c}^\transpose \vars^\star - \ubcost(\vars^\star)\\
              &\le \epsilon + \ubcost(\vars_\epsilon) - \tilde{c}^\transpose \vars_\epsilon + \tilde{c}^\transpose \vars^\star - \ubcost(\vars^\star) \\
              &\le \epsilon + \underset{\vars \in \varspace}{\max} \{ \ubcost(\vars) - \tilde{c}^\transpose \vars\} + \tilde{c}^\transpose \vars^\star- \ubcost(\vars^\star)\\
              &= \epsilon + \smallloss_{\GEN}(\tilde{c}).
        \end{align*}
        Here, the first inequality comes from the $\epsilon$-subptimality of $\vars_\epsilon$, the second inequality from the definition of the maximum over $\vars \in \varspace$.
        For item~\ref{thm:derivative:3.5}, if $\smallloss_{\GEN}(\predcost) = 0$, then for every $\vars \in \varspace$ by optimality of $\vars^\star$ in $\ubcost$, $
            \predcost^\transpose (\vars^\star - \vars) \le \ubcost(\vars^\star) - \ubcost(\vars) \le 0.$
        Therefore, $\predcost^\transpose (\vars - \vars^\star) \ge 0$ for every $\vars \in \varspace$.
        For uniqueness observe that the second inequality above holds strictly.
        Now we prove item~\ref{thm:derivative:5}. Fix $\epsilon > 0$. 
        By Lemma~\ref{lemma:sharpmin}, we can choose $\tilde{\truelayercost}$ such that $\vars^\star$ minimizes $\vars \mapsto \tilde{\truelayercost}^\transpose \vars$ over $\vars$ and such that for any $\vars \in \varspace,$ 
        $\tilde{\truelayercost}^\transpose (\vars - \vars^\star) \ge \ubcost(\vars) - \ubcost(\vars^\star) - \epsilon$. Therefore,
        $\smallloss_{\GEN}(\tilde{\truelayercost}|\ubcost) = \underset{\vars \in \varspace}{\max} \{ \ubcost(\vars) - \tilde{\truelayercost}^\transpose \vars \} + \tilde{\truelayercost}^\transpose \vars^\star - \ubcost(\vars^\star) \le \epsilon.
    $
        \ifarxiv
        \end{proof}
    \else
        \qedsymbol
    \endproof
\fi

We now discuss some choices for $\ubcost$ which lead to different loss functions with different interpretations and properties, which are summarized in Table~\ref{tab:losses}.
\begin{table}
    \centering
    \small
    \begin{tabular}{lllll}
    \toprule
        Loss function $\smallloss$
        & $\ubcost(\vars)$
        & $\smallloss = 0 \implies \smallloss_{\SUB} = 0$
        & MILP policy~\eqref{eq:costform1}
        & $\Omega(\vars)$ \\
        \midrule
        $\gspo$ & \eqref{eq:truecost} & \cmark & \cmark & 0\\
        $\ASL$ & $\| \vars - \vars^\star\|_2$ & \cmark & \cmark& 0\\
        $\text{Z}$ & 0 & If $\predcost \neq 0$ & \cmark & 0\\
        $\FY$ & 0 & \cmark & \xmark & $\|x\|_2^2$ \\
        \bottomrule
    \end{tabular}
    \label{tab:losses}
    \caption{A table summarizing properties of the loss functions discussed in Section~\ref{sec:surrogate}.}
\end{table}

    \subsection{Exact penalty}
    If $\ubcost(\vars) = \cost_\param(\vars)$ as in~\eqref{eq:truecost}, the loss function becomes
    \begin{equation} \label{eq:gspo} \ifarxiv \tag{\text{GSPO+}} \fi
    \smallloss_{\gspo}(\predcost) = \textstyle \underset{\vars \in \varspace}{\max} \{ \cost_\param(\vars) - \predcost^\transpose \vars\} + \predcost^\transpose \vars^\star- \cost_\param(\vars^\star),
    \end{equation}
    which we refer to as the \emph{generalized SPO+} (GSPO+) loss because it is inspired by the SPO+ loss~\cstmcite{spo} and it includes the nonlinear term $\cost_\param(\vars)$.
    In this case, solving the inner minimization problem in~\eqref{eq:truecost} is in general challenging.
    For this reason, we introduce the following approximations where the inner problem is easier to solve in practice. 
\subsection{Augmented suboptimality}
Take $\ubcost(\vars) = \nu d(\vars, \vars^\star)$ for a constant $\nu \in \reals_+$ and we recover the augmented suboptimality loss~($\ASL$)~\cstmcite{asl},
\begin{equation} \label{eq:asl} \ifarxiv \tag{$\text{ASL}$} \fi
    \smallloss^{\nu}_{\ASL}(\predcost) = \underset{\vars \in \varspace}{\max} \{\nu d(\vars, \vars^\star) - \predcost^\transpose \vars \} + \predcost^\transpose \vars^\star.
\end{equation}
Note that this version of the loss function has no dependence on the parameters~$\lowervars$ of the second-layer.
Optimizing this function over $\predcost$ can be interpreted as solving an inverse optimization problem where we try to make $\vars^\star$ optimal for the linear objective $\vars \mapsto \predcost^\transpose \vars$ without using knowledge of the second layer. 
\subsection{No penalty}
If $\ubcost(\vars) = 0$, loss function becomes
\begin{equation} \label{eq:oloss} \ifarxiv \tag{$\text{Z}$} \fi
    \smallloss_{\text{Z}}(\predcost) = \underset{\vars \in \varspace}{\max} \{- \predcost^\transpose \vars\} + \predcost^\transpose \vars^\star,
\end{equation}
which we refer to as the zero loss, \ie, $\text{Z}$-loss.
Clearly the $\text{Z}$-loss is minimized by $\predcost = 0$.
However it is also minimized by any $\predcost$ such that $\vars^\star \in \argmin \{\predcost^\transpose \vars \mid \vars \in \varspace\}$.
Therefore the effectiveness of a predictor that optimizes this loss greatly varies, depending on the local minimum it corresponds to.
    \subsection{Fenchel-young loss functions}
    While in this paper we focus on learning linear cost functions $\vars \mapsto \predcost^\transpose \vars$, which correspond to the MILP-based policies of the form~\eqref{eq:costform1}, the methods we propose can be extended to more generic policies of the form
    \begin{equation} \label{eq:nonlinearpredictor}
        \hat{\vars}(\param) = \begin{array}[t]{ll}
        \argmin & \predcost(\param)^\transpose \vars + \fenchel(\vars)
        \qquad  \text{subject to} \quad \vars \in \varspace(\param),
    \end{array}
    \end{equation}
    where~$\fenchel$ is a given convex function~\cstmcite{Dalle2022LearningWC}.
    To train the policy~\eqref{eq:nonlinearpredictor}, we need to modify our generic loss~\eqref{eq:linloss} by adding the convex penatly term $\fenchel$~\cstmcite{Dalle2022LearningWC},
    \begin{equation} \label{eq:genloss} \ifarxiv \tag{$\GEN$} \fi
        \smallloss^{\fenchel}_{\GEN}(\predcost) = \underset{\vars \in \varspace}{\max} \Big\{ \ubcost(\vars) - \ubcost(\vars^\star) +  \predcost^\transpose (\vars - \vars^\star) + \fenchel(\vars^\star) - \fenchel(\vars)\Big\}.
    \end{equation}
    As in~\eqref{eq:subloss} the optimization problem in~\eqref{eq:genloss} is nonconvex in $\vars$.
    However, as a maximum of convex functions, $\smallloss^{\fenchel}_{\GEN}(\predcost)$ is convex in $\predcost$.
    Note that the other loss functions discussed in this paper are of the form~\eqref{eq:genloss} with $\fenchel(\vars) = 0$ for every $\vars$.

    If we set $\ubcost = 0$, we recover the Fenchel-Young loss~\cstmcite{Blondel2019LearningWF} with penalty $\fenchel$,
    \begin{equation*}
        \smallloss^\fenchel_{\FY}(\predcost) = \underset{\vars \in \varspace}{\max} \big\{ \fenchel(\vars^\star) - \fenchel(\vars) +  \predcost^\transpose (\vars - \vars^\star) \big\}.
    \end{equation*}
    A convenient property of the Fenchel-Young loss, as opposed to our generic loss~\eqref{eq:linloss}, is that the maximization problem maximizes a convex objective which may simplify the computations at each training step.
    However, to evaluate the corresponding policy online, we need to solve~\eqref{eq:nonlinearpredictor}, which is often a mixed-integer convex optimization problem (MICP), as opposed to a MILP in~\eqref{eq:costform1}.
    This is why, we focus on MILP-based policies in this paper, but we still evaluate Fenchel-Young loss functions in our experiments.


\subsection{Effect of approximation error}
Theorem~\ref{thm:derivative} states that the inner optimization problem in~\eqref{eq:linloss} can be solved approximately to obtain an approximate gradient.
However, our generic loss~\eqref{eq:linloss} still depends on the exact optimizer $\vars^\star$ of $\ubcost$ and, consequently, of~\eqref{eq:master}.
Here we investigate the effect of $\vars^\star$ being misspecified.
Suppose that the $\vars^\star$ in~\eqref{eq:linloss} is replaced by $\vars^\star_\epsilon$ which is $\epsilon$-optimal for $\ubcost$. 
Let 
the misspecified general loss~\eqref{eq:linloss} be $\smallloss^\epsilon_{\GEN}(\predcost) = \underset{\vars \in \varspace}{\max} \{ \ubcost(\vars) - \ubcost(\vars^{\star}_\epsilon) - \predcost^\transpose (\vars - \vars^{\star}_\epsilon) \}$.

\begin{theorem} \label{thm:approxerror}
    Suppose $\vars^{\star}_{\epsilon}$ is $\epsilon$-optimal for $\ubcost$ and that $\predcost$ satisfies $\smallloss^\epsilon_{\GEN}(\predcost) \le \delta$.
    Then, for any $\vars \in \varspace^\star(\predcost)$, $\ubcost(\vars) \le \ubcost(\vars^{\star}) + \epsilon + \delta$. 
\end{theorem}
\ifarxiv
\begin{proof}
\else
\proof{Proof of Theorem~\ref{thm:approxerror}.}
\fi
The following is true,
\begin{equation*}
    \delta \ge \smallloss^\epsilon_{\GEN}(\predcost)  
    \ge \ubcost(\vars) - \predcost^\transpose (\vars - \vars^{\star}_{\epsilon}) - \ubcost(\vars^{\star}_{\epsilon})
    \ge \ubcost(\vars) - \ubcost(\vars^\star) - \epsilon.
\end{equation*}
The first inequality is the assumption. The second is by definition of the maximum. The third is because $\vars \in \varspace^\star(\predcost)$ and $\vars^{\star}_{\epsilon}$ is $\epsilon$-optimal for $\ubcost$.
\ifarxiv
\end{proof}
\else
\qedsymbol
\endproof
\fi
This result means that we do not need to solve problem~\eqref{eq:master} to optimality. Instead, we can bound the prediction error as the sum of the suboptimality $\epsilon$ of the approximate solution and the value of the loss function $\delta$.

\section{Learning with suboptimality bounds using conformal prediction}
\label{sec:lowerbounds}
We assume to have access to a \emph{training} dataset $\dataset = \{\param_i, (\vars^{\star}_i, \lowervars^\star_i)\}_{i =1}^N$ where each $\param_i$ is an iid draw from distribution $\paramdist$ and each pair $(\vars^\star_i, \lowervars^\star_i)$ solves problem~\eqref{eq:master} given~$\theta_i$.
We train a predictor $\nn_w$ with learned parameters $w$ by minimizing an empirical approximation of the suoptimality loss~\eqref{eq:linloss}.
However, as discussed in Section~\ref{sec:surrogate}, loss function $\smallloss_{\SUB}$ may be nonconvex and discontinuous and we approximate it with a convex surrogate function.
Specifically, we train by minimizing the following function using stochastic gradient descent,
\begin{equation} \label{eq:approxempiricalrisk}
    \hat{\risk}(\nn_w) = \frac{1}{N} \sum_{i = 1}^N \smallloss(\nn_w(\param_i)|\param_i),
\end{equation}
where $\smallloss$ is one of the approximate loss functions in Table~\ref{tab:losses}.
By this method, we can learn an optimizer which provides heuristic feasible solutions and upper bounds for our problem. 

By solving~\eqref{eq:costform1} and~\eqref{eq:lowerprob}, our trained predictors always provide feasible solutions~$(\hat{\vars}, \hat{\lowervars})$ for the original problem~\eqref{eq:master}, and therefore valid upper bounds on its optimal objective.
However, our architecture does not provide valid lower bounds on the true cost, which are necessary to quantify the suboptimality of the predicted solutions.
In this section, we provide online suboptimality bounds using conformal prediction~\cstmcite{conformal}.
\subsection{Conformal prediction-based bounds}
\label{sec:conformal}
Assume that we have access to a \emph{calibration} dataset $\calibration = \{\param_i, (\vars^{\star}_i, \lowervars^\star_i)\}_{i = 1}^M$ that is independent of $\dataset$, where
$\param_i$ are iid according to $\paramdist$ and $(\vars_i^\star, \lowervars_i^\star)$ solve~\eqref{eq:master} given $\param_i$.

Suppose that we have trained the predictor $\predcost_\nnparam(\param)$ on the training dataset $\dataset$.
Given $\param$, let $\ub(\param)$ be the upper bound on the true cost $z(\param)$ from~\eqref{eq:master} given by the feasible solution $(\hat{\vars}, \hat{\lowervars})$ outputted by our method, \ie, $\ub(\param) = \costx(\param)^\transpose \hat{\vars}(\param) + \costy(\param)^\transpose \hat{\lowervars}(\param)$.
Let $\lb(\param)$ be any function of $\param \in \paramspace$ which gives a lower-bound on the true cost $z(\param)$.
In the case of MILPs, for example, this function could be the optimal value of the continuous relaxation of~\eqref{eq:master}.
Assume that $\ub$ and $\lb$ are independent of $\calibration$.

For any $\ub > \lb \in \reals$, let $\phi_{l,u} : [\lb, \ub) \rightarrow \reals_+$ be a monotonic increasing homomorphism.
Let $\conff$ be any non-negative function of $\param \in \paramspace$ which is independent of the calibration dataset $\calibration$ and such that $\conff(\param) \in [\lb(\param), \ub(\param))$ for every $\param$.
In our experiments we take $\phi_{l,u}(x) = \operatorname{arctanh}\left((x - l) / (u - l)\right)$.
\begin{theorem} \label{thm:conformal}
Fix $\alpha \in (0, 1)$.
For all samples $i=1,\dots,M$ in the calibration dataset $\calibration$, let
\begin{equation*}
    \Phi_i = \frac{\phi_{\lb(\param_i), \ub(\param_i)}\big(\conff(\param_i)\big)}{\phi_{\lb(\param_i), \ub(\param_i)}\big(z(\param_i)\big)}.
\end{equation*}
Let $\{\Phi_{[1]}, \dots \Phi_{[M]}\}$ denote the values of~$\Phi_i$ sorted in increasing order, so that $\Phi_{[1]} \le \Phi_{[2]} \le \dots \le \Phi_{[M]}$.
Define
$\nu_\alpha = \lceil M\alpha \rceil/(M+1)$
and let $q_\alpha = \Phi_{\nu_\alpha (M+1)}$.
Then, 
\begin{equation*}
\prob\left[\frac{\phi_{\lb(\param), \ub(\param)}( \conff(\param)\big)}{\phi_{\lb(\param), \ub(\param)}(z(\param))}  \le q_\alpha  \right] = 1 - \frac{\lceil M\alpha \rceil}{M+1},
\end{equation*}
where the probability is taken over randomness in the calibration dataset $\calibration$ and a test parameter $\param \sim \paramdist$ that is independent of $\calibration$.
\end{theorem}
This theorem is a direct application of conformal prediction guarantees~\cstmcite{conformal,conformal_vovk}.
It is based on the following lemma about the ordering of random variables.
\begin{lemma} \label{lemma:rvorder}
    Let $\alpha \in (0,1)$ and $X_1, \dots X_n, X_{n+1}$ be iid continuous random variables taking values in $[0, \infty]$. Let $Q$ be the random variable given by the $\lceil n\alpha \rceil / n $ quantile of $\{X_1, \dots X_n\}$. Then,
    $
        \mathbf{P}(X_{n+1} \le Q) = 1 - \lceil n\alpha \rceil / (n+1).
    $
\end{lemma}
\ifarxiv
\begin{proof}
\else
\proof{Proof of Lemma~\ref{lemma:rvorder}}
\fi
    Since $X_1, \dots X_{n+1}$ are iid and continuous, all orderings of $X_1, \dots X_{n}$ are equally likely, and the probability that any pair is equal is 0. Given~$q \in \{0, \dots n\}$, let $E$ be the event $E = \{X_{n+1} \text{ is less than at least $n - q$ of }X_1, \dots X_n \}$.
    The probability that $X_{n+1}$ is in position $j$ in the ordering of all the $X_i$ is exactly $1/(n+1)$ for each $j \in \{1, \dots n+1\}$. Therefore, the probability that it is in the bottom $q$ is exactly $q/(n+1)$ and so $\mathbf{P}(E) = q/(n+1)$.
    The statement of the lemma is exactly this with $q = \lceil n \alpha \rceil.$
\ifarxiv
\end{proof}
\else
\qedsymbol
\endproof
\fi
We can now prove Theorem~\ref{thm:conformal}.
\ifarxiv
\begin{proof}
\else
\proof{Proof of Theorem~\ref{thm:conformal}.}
\fi
    This is a consequence of Lemma~\ref{lemma:rvorder} with $n=M$ and,
    \begin{equation*}
        X_i = \Phi_i, \quad 
        X_{n+1} = \frac{\phi_{l(\param), u(\param)}\big(h(\param)\big)}{\phi_{l(\param), u(\param)}\big(z(\param)\big)},
    \end{equation*}
    since the random variable $Q$ is given by $q_\alpha$.
    If any $\phi_{l(\param_i), u(\param_i)}\big(z(\param_i)\big) = 0$ set $X_i = \infty$.
\ifarxiv
\end{proof}
\else
\qedsymbol
\endproof
\fi
This theorem provides a systematic way to construct bounds on the true optimal value $z(\param)$ for a given parameter $\param \sim \paramdist$.
We now clarify the specific choices of the functions $\conff$, $\lb$, $\ub$, and $\phi$.
Although the result holds for arbitrary functions, choosing $h(\param_i)$ close to $z(\param_i)$ yields tighter bounds.
Our goal is to select functions that incorporate information from the learned predictor $\predcost$ to produce informative bounds on the optimality gap of $\hat{\vars}(\param)$.
We remark that the functions $\conff$, $\lb$, and $\ub$ may depend on the training dataset $\dataset$ or any other dataset independent of $\calibration$ without violating the assumptions of Theorem~\ref{thm:conformal}.

\subsubsection{Training the conformal predictor}
We parametrize the function $\conff$ by a neural network $\confnn_\confparam$ with weights $\confparam$.
To train this network, we introduce a third dataset called the \emph{evaluation dataset}, denoted $\eval = \{\param_i, (\vars_i^\star, \lowervars_i^\star)\}_{i = 1}^E$, where $\param_i$ are iid draws from $\paramdist$ and $(\vars_i^\star, \lowervars_i^\star)$ solve~\eqref{eq:master} for $\param_i$.
The dataset $\eval$ is independent of both the training dataset $\dataset$ and the calibration dataset $\calibration$.
The neural network $\confnn_\confparam$ takes as input the parameter $\param \in \paramspace$, the predicted solution $(\predvars(\param), \predlowervars(\param))$, and the available bounds $\lb(\param)$ and $\ub(\param)$.
It outputs a value $\conff(\param) = \confnn_\confparam(\param, \hat{\vars}(\param), \hat{\lowervars}(\param), \lb(\param), \ub(\param))$ that is constrained to lie in $[\lb(\param), \ub(\param)]$ via a sigmoid activation in the final layer.
We train the weights $\confparam$ to predict the true optimal value $z(\param)$ by minimizing the squared error over the evaluation dataset,
\begin{equation} \label{eq:confloss} \ifarxiv \tag{$\hat{\risk}_{\text{CONF}}$} \fi
    \frac{1}{E} \sum_{i = 1}^E \left(z(\param_i) -  \confnn_\confparam\big(\param_i, \hat{\vars}(\param_i), \hat{\lowervars}(\param_i), \lb(\param_i), \ub(\param_i)\big)\right)^2.
\end{equation}
Once $\confnn_\confparam$ is trained, we use the calibration dataset $\calibration$ to obtain a threshold $q_\alpha$ for a desired confidence level~$\alpha$ as described in Theorem~\ref{thm:conformal}.
This threshold enables us to compute online bounds as explained in the following section.

\subsubsection{Online bounds in probability}
Let $\predcost_\nnparam(\param)$ be a learned function trained to minimize the risk~\eqref{eq:approxempiricalrisk} for some approximate loss function from Table~\ref{tab:losses}.
Given the predictor $\predcost_\nnparam(\param)$, let $\hat{\vars}(\param)$ denote the predicted optimal solution to~\eqref{eq:costform1}.
Let the conformal predictor $\confnn_\confparam$ be trained to minimize~\eqref{eq:confloss} using $\predcost_\nnparam(\param)$ and bounding functions $\ub(\param)$ and $\lb(\param)$.
For a given confidence level $\alpha \in (0,1)$, let $q_\alpha$ be the threshold obtained from Theorem~\ref{thm:conformal}.
Algorithm~\ref{alg:hierarchical-mip-solution-conformal} describes a procedure to efficiently compute a high-probability bound on the suboptimality of $\hat{\vars}$.
\begin{algorithm}
\caption{Hierarchical MIP Solution with Conformal Bound} \label{alg:hierarchical-mip-solution-conformal}
\textbf{Input:} Parameter $\param$, predictor $\predcost_\nnparam(\param)$, conformal predictor $\confnn_\confparam(\param)$, lower-bound function $\lb(\param)$, upper-bound function $\ub(\param)$, quantile $q_\alpha$ from calibration.
\begin{algorithmic}[1]
    \State $\predvars(\param) \gets$ solve upper-level problem~\eqref{eq:costform1}
    \State $\predlowervars(\param) \gets$ solve lower-level problem~\eqref{eq:lowerprob} with $\vars = \hat{\vars}(\param)$
    \State $\lb(\param) \gets$ solve convex relaxation of~\eqref{eq:master}\Comment{lower bound}
    \State $\ub(\param) \gets \costx(\param)^\transpose \hat{\vars}(\param) + \costy(\param)^\transpose \hat{\lowervars}(\param)$ \Comment{upper bound from feasible solution}
    \State $\conff(\param) \gets \confnn_\confparam(\param, \hat{\vars}(\param), \hat{\lowervars}(\param), \lb(\param), \ub(\param))$ \Comment{predicted optimal solution value}
    \State $\omega(\param) \gets \phi^{-1}_{\lb(\param), \ub(\param)} \left( q \cdot \phi_{\lb(\param), \ub(\param)}(\conff(\param)) \right)$ \Comment{predicted bound}
\end{algorithmic}
\textbf{Output:} Feasible solution $(\hat{\vars}(\param), \hat{\lowervars}(\param))$ and probabilistic lower bound $\omega(\param)$. 
\end{algorithm}
This procedure yields a feasible solution $\hat{\vars}(\param)$ to the upper-level problem~\eqref{eq:costform1} and a bound $\omega(\param)$ such that with probability at least $1 - \lceil M \alpha \rceil/(M+1)$, the true optimal value satisfies $z(\param) = \costx(\param)^T \vars^\star(\param) + \costy(\param)^T \lowervars^\star(\param) \ge \omega(\param)$.
In addition, its computational cost is much lower than solving the full problem~\eqref{eq:master}.
\begin{figure}
    \centering
    \includegraphics[trim={0 10cm 0 0cm},clip,width=0.8\linewidth]{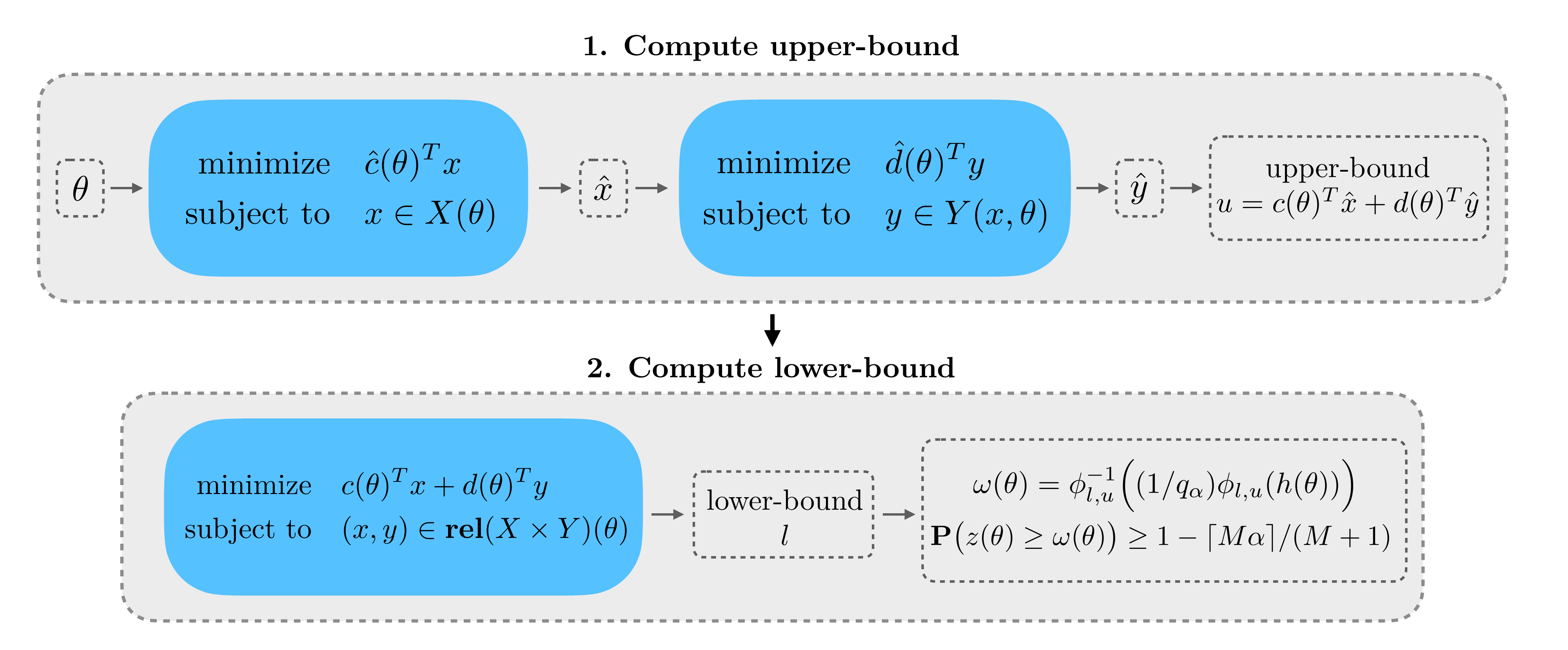}
    \caption{The full architecture of our model to predict a feasible solution to the top problem $\hat{\vars}$ and obtain a bound for its suboptimality.
    Here we write $\textbf{rel}(X)$ to mean the convex set given by relaxing integrality constraints in the set $X$.}
    \label{fig:full}
\end{figure}
\section{Computational experiments} \label{sec:computation}
We demonstrate the effectiveness of our methods on a set of standard baselines.

\paragraph{Hardware.} 
All experiments were run on one of Princeton OIT's computer clusters, which consists of a range of Intel, AMD and ARM processors.
The experiments were run on CPUs and each experimental run used no more than 16GB RAM.
All experiments are single-threaded. 

\paragraph{Data.}
For each experiment, we generate $\NINSTANCE$ problem instances and solve them to global optimality.
We partition this data into four subsets: a training dataset $\dataset$ with $N=\NTRAIN$ samples, an evaluation dataset $\eval$ with $E=\NEVAL$ samples, a calibration dataset $\calibration$ with $M=\NCAL$ samples, and a test dataset $\test = \{\param_i, (\vars_i^\star, \lowervars_i^\star)\}_{i = 1}^T$ with $T=\NTEST$ samples.
The training dataset $\dataset$ is used to train the cost predictor $\predcost_\nnparam(\param)$ by minimizing the loss functions in Table~\ref{tab:losses}.
The evaluation dataset $\eval$ serves two purposes: first, to perform grid search over training hyperparameters such as the learning rate (details are available in the shared repository), and second, to train the conformal predictor $\confnn_\confparam$ by minimizing the conformal loss~\eqref{eq:confloss}.
The calibration dataset $\calibration$ is used to calibrate the conformal predictor as described in Section~\ref{sec:lowerbounds}.
We report performance metrics on the test set $\test$, which remains unseen by all methods during training.
All problem instances are randomly generated with parameter dimension $p=\NTHETA$, but the optimization problem sizes vary across experiments.

\paragraph{Evaluation.}
To evaluate feasible solution quality, we consider the regret at a feasible solution $\hat{\vars}(\param)$ given by $\cost_\param(\hat{\vars}(\param)) - \cost_\param(\vars^\star(\param))$ with $f_\theta$ defined in~\eqref{eq:truecost}.
We compare methods based on their average regret over the test dataset~$\test$, 
\begin{equation} \label{eq:absgap}
    \hat{R}_{\text{ABS}} = \frac{1}{T} \sum_{i = 1}^T \cost_{\param_i}(\hat{\vars}) - \cost_{\param_i}(\vars^\star(\param_i)),
\end{equation}
and the normalized regret given by,
\begin{equation*}
    \hat{R}_{\text{NORM}} = \frac{1}{T} \sum_{i = 1}^T \frac{\cost_{\param_i}(\hat{\vars}) - \cost_{\param_i}(\vars^\star(\param_i))}{\cost_{\param_i}(\vars^\star(\param_i))}. 
\end{equation*}
We compare the time required to compute a feasible solution $(\hat{\vars}, \hat{\lowervars})$ for~\eqref{eq:master} across different baselines.
For our method, we first solve~\eqref{eq:costform1} to obtain $\predvars$, then sequentially solve each lower-level problem~\eqref{eq:lowerprob} to obtain $\predlowervars$.
Note that the lower-level problems can be solved in parallel, which would significantly reduce computation time, but we do not exploit parallelization in our experiments.
To evaluate the quality of the predicted lower bounds $\omega(\param)$ produced by our conformal predictor, we use the following metrics,
\begin{align*}
        \hat{r}_{\text{rel}}^+ &= \frac{1}{T} \sum_{i =1}^T \left(\frac{z(\param_i) - \omega(\param_i)}{z(\param_i) - \omega_{\text{rel}}(\param_i)}\right) \mathbf{1}_{z(\param_i) - \omega(\param_i) \ge 0}, \\
        \hat{r}_{\text{rel}}^- &= \frac{1}{T} \sum_{i =1}^T \left(\frac{-z(\param_i) + \omega(\param_i)}{z(\param_i) - \omega_{\text{rel}}(\param_i)}\right) \mathbf{1}_{z(\param_i) - \omega(\param_i) \le 0}, \\
        \hat{r}_{\%} &= \frac{1}{T} \sum_{i =1}^T \mathbf{1}_{z(\param_i) - \omega(\param_i) \le 0},
\end{align*}
where $\omega_{\text{rel}}(\param)$ is the objective value of the convex relaxation and $\mathbf{1}$ denotes the indicator function.
Here $\hat{r}_{\text{rel}}^+$ measures the average relative optimality gap when the bounds are valid, $\hat{r}_{\text{rel}}^-$ measures the average relative violation when the bounds are invalid, and $\hat{r}_{\%}$ is the fraction of invalid bounds.
We plot the evaluation regret over time for each learned method and display the time required to solve all test instances.
We present violin plots of the bound quality from conformal prediction, training the conformal predictor only on the method with the best performance on the evaluation set.
All numerical results are also provided in tables.
We use $\alpha=\CPALPHA$ for the conformal prediction experiments.

\paragraph{Baselines.}
We compare the solution quality of our method against a range of global solvers and heuristics.
The following methods are evaluated on each problem instance.
\begin{itemize}
    \item \GRB: The Gurobi~\cstmcite{gurobi} solver, configured to solve each instance to optimality within a tolerance of $\TOLERANCE$ and a time limit of $\TIMELIMIT$ seconds.
    \item \GRBo: The Gurobi solver configured to terminate after finding a single feasible solution. In plots, we also show a version that terminates after finding three feasible solutions.
    \item SCIP: The \SCIP~\cstmcite{Achterberg2009SCIPSC} solver, configured to solve each instance to optimality within a tolerance of \TOLERANCE~and a time limit of $\TIMELIMIT$ seconds.
    \item \SCIPo: The SCIP solver configured to terminate after finding a single feasible solution. In plots, we also show a version that terminates after finding three feasible solutions.
    \item Nearest neighbor (NN): For a new parameter $\param$, this method finds the index $i \in \{1, \dots, N\}$ that minimizes $\|\param - \param_i\|_2^2$ over the training dataset $\dataset$, then returns the Euclidean projection of $\vars^\star_i$ onto the feasible region of~\eqref{eq:costform1}.
    \item Direct prediction (DP): We train a neural network $r_w(\param)$ with weights $w$ to predict $\vars^\star(\param)$ directly by minimizing the squared loss $(1/N)\sum_{i=1}^N \|r_w(\param_i) - \vars^\star(\param_i)\|_2^2$ over the training dataset.
    For a new parameter $\param$, the method returns the Euclidean projection of $r_w(\param)$ onto the feasible region of~\eqref{eq:costform1}.
    \item Our approach: We train the cost predictor $\predcost(\param)$ by minimizing the $\ASL$ and $\text{Z}$-losses from Table~\ref{tab:losses}.
    We also train a predictor using the FY loss with convex penalty $\fenchel(\vars) = \|\vars\|_2^2$.
    At test time, we use Gurobi to first solve~\eqref{eq:costform1} to obtain $\predvars$ (or~\eqref{eq:nonlinearpredictor} for the FY loss), then solve~\eqref{eq:lowerprob} to obtain $\predlowervars$, as illustrated in Figure~\ref{fig:full}.
    When the lower-level problems are separable, we solve them sequentially. The reported time to find a feasible solution includes the time to solve all lower-level problems.
\end{itemize}

\paragraph{Learned models.} 
We parametrize the cost vector $\predcost_\nnparam(\param)$ for $\vars$ as a neural network that takes $\param$ as input.
We parametrize the conformal predictor $\confnn_\confparam$ as a neural network that takes as input $\param$, $\predvars$, $\predlowervars$, $l$, and $u$.
Both neural networks are $\NLAYER$-layer feedforward networks with ReLU activation functions and $\NNEURON$ neurons per layer.
The conformal predictor $\confnn_\confparam$ uses a sigmoid activation function after the final layer to ensure its output lies in $[l, u]$, and we take $\phi_{l,u}(x) = \operatorname{arctanh}\left((x - l) / (u - l)\right)$.

\paragraph{Code.}
The code for our experiments is available at
\begin{center}
    \url{anonymous.4open.science/r/hmip-23B9}.
\end{center}
We first discuss the results for the experiments relating to learning the hierarchical models and then present results for the conformal prediction-based bounds.

\subsection{Hierarchical knapsack problems}
We consider the following hierarchical knapsack problem,
\begin{equation*}  \ifarxiv \tag{$\mathcal{P}_\text{K}$} \fi
    \begin{array}{ll}
        \text{minimize} \quad & \sum_{j=1}^J \costy_j(\param)^\transpose \lowervars_j \\
        \text{subject to} \quad & a_j^\transpose \lowervars_j \le b_j, \quad j=1, \dots J, \\
        & \lowervars_j \le \vars_j \ones, \quad j=1, \dots J, \\
        & a_0^\transpose \vars \le b_0, \quad \vars \in \{0, 1\}^J, \quad \lowervars_j \in \{0,1\}^{k}, \quad j=1, \dots J.
    \end{array}
\end{equation*}
There are $J \in \integers_+$ lower-level knapsacks, each of size $k \in \integers_+$, with capacity $b_j \in \reals_+$ and weight vector $a_j \in \reals_+^k$ for $j=1,\dots,J$.
The vector $\ones$ denotes the vector of all ones in $\reals^k$.
The cost vector for lower knapsack $j$ is given by $\costy_j(\param) \in \reals^{k}$ and depends on the parameter $\param$.
There is a single upper-level knapsack with weights $a_0 \in \reals_+^J$ and capacity $b_0 \in \reals_+$.
For each $j \in \{1, \dots, J\}$, lower knapsack $j$ may only be filled if item $j$ is selected in the upper knapsack.
The optimization variables are $\vars \in \{0,1\}^J$ and $\lowervars_1, \dots, \lowervars_J \in \{0,1\}^k$.
We seek to learn a cost vector $\predcost_\nnparam(\param)$ that predicts the optimal solution for the upper knapsack given $\param$ without solving for the lower knapsack variables.
This reduces the number of variables from $(k+1)J$ to $J$.
Once we obtain a feasible solution for the upper knapsack, the lower knapsacks become separable independent problems as in~\eqref{eq:decomp}.
The upper-level problem is therefore
\begin{equation*}  \ifarxiv \tag{$\hat{\mathcal{P}}_\text{K}$} \fi
    \begin{array}{ll}
        \text{maximize} \quad & \predcost(\param)^\transpose \vars\\
        \text{subject to} \quad
        & a_0^\transpose \vars \le b_0, \quad
        \vars \in \{0, 1\}^J.
    \end{array}
\end{equation*}

\paragraph{Problem generation.}
We generate $b_0, \dots b_h$ and $a_0, \dots a_h$ as the absolute values of standard gaussian distributions.
We also generate $A \in \reals^{Jk \times p}$ according to a standard normal distribution. 
To create the parametric family we generate $\NINSTANCE$ values of $\param \in \reals^p$ according to the absolute value of a standard normal distribution and let $d(\param) = -|A \param| \in \reals^{Jk}$. We take $J = k = \NKNAPSACK$.

\begin{figure}
     \centering
     \includegraphics[scale=0.8, trim={0 0.3cm 0 0}]{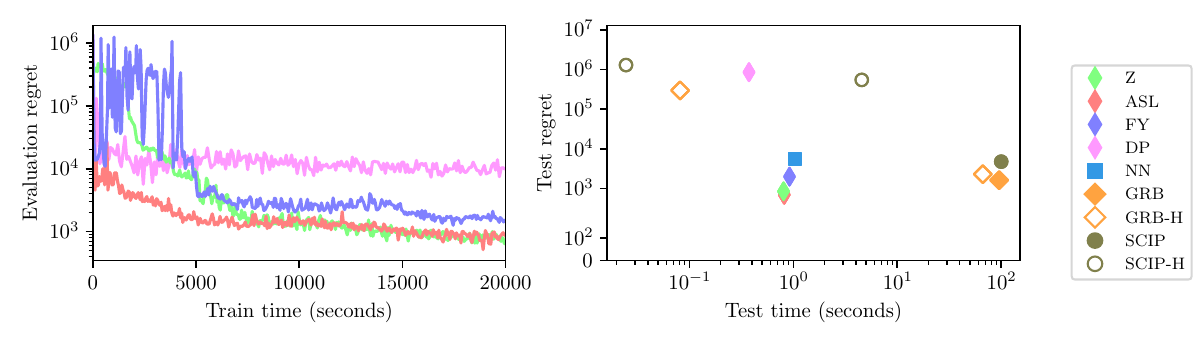}
     \caption{Results for the knapsack experiment. The left plot shows evaluation regret throughout the training process. The right plot shows the test regret~\eqref{eq:absgap} and time (averaged over the test instances) on the test set.
     The first heuristic to terimnate for Gurobi and SCIP terminates after finding a single feasible solution. The second to terminate terminates after finding three.} 
    \label{fig:knapsack}
\end{figure}

\paragraph{Results.} 
Training loss, evaluation regret, and test dataset results are shown in Figure~\ref{fig:knapsack}.
Our method finds high-quality feasible solutions significantly faster than Gurobi and SCIP.
The solution obtained after one second is on average better than the solution found by Gurobi after nearly 100 seconds.
The feasible solutions are also higher quality than those found by the nearest-neighbor method.

\begin{table}[ht]
\centering
{\fontsize{7}{8.5}\selectfont
\begin{tabular}{l r r r r r r r r r}
\toprule
 & NN & DP & GRB & GRB-H & SCIP & SCIP-H & \textbf{ASL} & \textbf{Z} & \textbf{FY} \\
 \midrule
  \begin{filecontents*}{knapsack_scientific.csv}
name,NN,DP,GRB,GRB_H,SCIP,SCIP_H,ASL,Z,FY
Time(s),1.04e+00,3.70e-01,9.55e+01,6.65e+01,1e+02,4.55e+00,8.20e-01,8.10e-01,9.20e-01
Regret,5.5e+03,8.5e+05,1.6e+03,2.3e+03,4.7e+03,5.4e+05,7e+02,8.5e+02,2e+03
Normalized Regret,4.40e-03,6.63e-01,1.30e-03,1.80e-03,3.70e-03,4.22e-01,6.00e-04,7.00e-04,1.90e-03
\end{filecontents*}
\csvreader[
      head to column names,
      late after line=\\,
      before reading={\catcode`\ =9\relax}
      ]{knapsack_scientific.csv}{}{%
      \csvcoli & \csvcolii & \csvcoliii & \csvcoliv & \csvcolv &
      \csvcolvi & \csvcolvii & \csvcolviii & \csvcolix & \csvcolx
      }
    \bottomrule
\end{tabular}
\caption{Results for the knapsack experiment. 
    The left plot shows evaluation regret throughout the training process.
     The right plot shows the test regret~\eqref{eq:absgap} and time (averaged over the test instances) on the test set.
     The first heuristic to terimnate for Gurobi and SCIP terminates after finding a single feasible solution. The second to terminate terminates after finding three.}
}
\end{table}

\subsection{Capacitated facility location problems}
We consider a family of capacitated facility location problems~\cstmcite{Cornujols1991ACO, conf/nips/GasseCFCL19}.
In this problem, $I$ is a set of clients and $J$ is a set of sites where facilities can be located.
The objective is to match all clients to facilities in a way that minimizes the total cost.
The variable $\vars_j \in \{0,1\}$ indicates whether the facility at location $j \in J$ is built, and the variable $\lowervars_{i,j} \in [0,1]$ represents the proportion of client $i$'s demand serviced by facility $j$.
The cost of building a facility at location $j$ is $\costx_j(\param) \in \reals_+$, and the cost of servicing customer $i$ from location $j$ is $\costy_{ij}(\param) \in \reals_+$.
We include a penalty $\gamma \sum_{i \in I} \faclocz_i$ for some fixed $\gamma > 0$, where $\faclocz_i$ represents the unmet demand for customer $i$.
The demand for customer $i$ is given by $e_{i}(\param) \in \reals_+$, and the capacity of location $j$ is given by $s_j \in \reals_+$.
The matrix $A \in \reals^{\NFACLOCCOMPLICATING \times |J|}$ and vector $b \in \reals^{\NFACLOCCOMPLICATING}$ represent complicating constraints.
The problem is formulated as
\begin{equation} \label{eq:facloc} \ifarxiv \tag{$\mathcal{P}_\text{F}$} \fi
    \begin{array}{ll}
        \text{minimize} & \sum_{i \in I} \sum_{j \in J} \costy_{ij}(\param) \lowervars_{ij} + \sum_{j \in J} \costx_j(\param) \vars_j + \gamma \sum_{i \in I} \faclocz_i \\
        \text{subject to} & \sum_{j \in J} \lowervars_{ij} = 1 - \faclocz_i, \quad  i \in I, \\
        & A \vars \le b, \\
        & \sum_{i \in I} e_i(\param) \lowervars_{i j} \le s_j \vars_j, \quad  j \in J, \\
        & \lowervars_{ij} \in [0, 1], \quad \vars_j \in \{0, 1\}, \quad \faclocz_i \in [0,1], \quad i \in I, j \in J.
    \end{array}
\end{equation}
We seek to predict which facilities to build without solving the entire problem~\eqref{eq:facloc}.
To form the upper-level problem, we drop the lower-level variables $\lowervars$ and penalties $\faclocz$ and introduce a predicted cost $\predcost(\theta)$ to obtain
\begin{equation} \label{eq:facloctop} \ifarxiv \tag{$\hat{\mathcal{P}}_\text{F}$} \fi
    \begin{array}{ll}
        \text{minimize} & \predcost(\param)^T \vars \\
        \text{subject to} & A \vars \le b, \quad \vars \in \{0,1\}^{|J|}.
    \end{array}
\end{equation}
Once we have solved~\eqref{eq:facloctop}, we obtain the complete solution by solving~\eqref{eq:facloc} with $\vars$ fixed, which corresponds to a continuous lower-level problem.

\paragraph{Problem generation.}
We set $|I|=|J|=\FACLOCN$ and generate problem instances according to the following procedure.
We draw the $\NINSTANCE$ parameter values $\param_i$ from a standard normal distribution.
We generate the matrix $A$ from the absolute value of a Gaussian distribution and set $b$ to be half the row sums of $A$.
We generate the parametric family and remaining parameters following~\cstmciter{conf/nips/GasseCFCL19}, but modify the procedure so that $\costy_{ij}(\param)$, $\costx_j(\param)$, and $e_i(\param)$ are approximately linear functions of $\param$.
We set $\gamma = \FACLOCGAMMA$.

\begin{figure}
    \centering
    \includegraphics[scale=0.8, trim={0 0.3cm 0 0}]{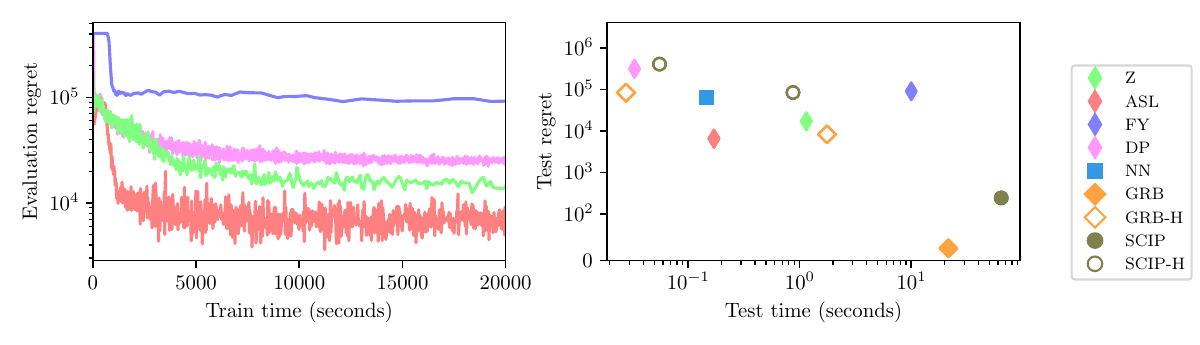}
    \caption{Results for the facility location experiment. 
    The left plot shows evaluation regret throughout the training process.
     The right plot shows the test regret~\eqref{eq:absgap} and time (averaged over the test instances) on the test set.
     The first heuristic to terimnate for Gurobi and SCIP terminates after finding a single feasible solution. The second to terminate terminates after finding three.}
    \label{fig:facloc}
\end{figure}

\paragraph{Results.}
Results are displayed in Figure~\ref{fig:facloc}.
The Z and ASL losses significantly outperform the nearest-neighbor method.
The ASL loss finds a feasible solution in approximately 0.1 seconds that is, on average, higher quality than the solution Gurobi obtains after one full second.
While SCIP and Gurobi heristics find feasible solutions quickly on average in the test set, these solutions are of low quality.
\begin{table}[ht]
\centering
{\fontsize{7}{8.5}\selectfont
\begin{tabular}{l*{9}{r}}
\toprule
 & NN & DP & GRB & GRB-H & SCIP & SCIP-H & \textbf{ASL} & \textbf{Z} & \textbf{FY}\\
\midrule
  \begin{filecontents*}{facloc_scientific_2.csv}
,NN,DP,GRB,GRB_H,SCIP,SCIP_H,ASL,Z,FY
Time,1.50e-01,3.00e-02,2.16e+01,1.76e+00,6.42e+01,8.70e-01,1.70e-01,1.15e+00,1.00e+01
Regret,6.47e+04,3.12e+05,2.63e+01,8.29e+03,2.46e+02,8.40e+04,6.55e+03,1.71e+04,9.01e+04
Normalized Regret,6.43e+04,4.49e+05,1.44e+00,6.97e+03,1.46e+02,9.34e+04,1.02e+04,1.86e+04,1.15e+05
\end{filecontents*}
\csvreader[
      head to column names,
      late after line=\\,
      before reading={\catcode`\ =9\relax}
      ]{facloc_scientific_2.csv}{}{%
      \csvcoli & \csvcolii & \csvcoliii & \csvcoliv & \csvcolv &
      \csvcolvi & \csvcolvii & \csvcolviii & \csvcolix & \csvcolx
      }
    \bottomrule
\end{tabular}
}
\caption{Results for the facility location experiment.
Our method names are in bold.
The heuristic methods recorded here terminate after finding a single feasible solution.}
\end{table}

\subsection{Multi-agent heterogeneous routing problems}
We consider a capacitated vehicle routing problem with heterogeneous vehicles and fuel constraints.
Let $V$ be a set of nodes and let $E$ be a set of directed edges between the nodes forming a directed graph $G=(V,E)$.
Let $v_0 \in V$ denote the depot, and let each vertex $v$ have an associated demand $\demand_v \in \reals_+$.
Suppose we have $m \in \integers_+$ vehicles, where vehicle $k \in \{1, \dots, m\}$ has capacity $\gamma_k \in \reals_+$ and starting fuel $f_k \in \reals_+$.
Traversing edge $e \in E$ requires $l_e \in \reals_+$ units of fuel.
The capacitated vehicle routing problem on graph $G$ is formulated as
    \begin{subequations}
        \label{eq:vrp} 
        \begin{align}
            \text{maximize} \quad & \textstyle  \sum_{v \in V} \demand_v(\param) \routingz_v \\
        \text{subject to} \quad 
          & \textstyle \sum_{k=1}^K \vars^{k}_{v} \le 1, \quad v \in V\setminus \{v_0\}, \label{eq:visitnodeonce} \\
        & \textstyle \sum_{v \in V} \demand_v(\param) \vars^k_{v} \le \gamma_k, \quad k=1, \dots m, \label{eq:capacity}\\
        & \textstyle \vars^k_{v} \le \sum_{e \in E^+_v} \routingw^k_{e}, \quad v \in V, \quad k=1, \dots m, \label{eq:vrp:3}\\
        & \textstyle \vars^k_{v_0} = 1, \quad k=1, \dots m, \label{eq:vrp:start}\\
        & A \routingw^{k} = 0, \quad k=1, \dots m, \label{eq:vrp:flow} \\
        & \textstyle \vars^k_{v} \le \sum_{u \in E^+_S}  \routingw^k_{u}, \quad k=1, \dots m, S \subseteq V, v \in S, v_0 \notin S \label{eq:vrp:subtour} \\
        & \textstyle \sum_{e \in E} l_e \routingw^k_{e} \le f_k, \quad  e \in E, \label{eq:vrp:fuel} \\
        & \textstyle \routingz_v \le \sum_{k=1}^m \vars^k_{v}, \quad  v \in V, \label{eq:vrp:reward}\\
        & \textstyle \routingz \in \{0, 1\}^V, \quad \routingw^k \in \{0, 1\}^E, \quad \vars^k \in \{0, 1\}^V,
    \end{align}
\end{subequations}
where $A \in \reals^{|V| \times |E|}$ is the directed incidence matrix of graph $G$, $E^+_v$ denotes the set of edges leaving node $v$, $\demand_v(\param) \in \reals_+$ is the demand at node $v$, $S$ represents a subset of nodes, and $E^+_S$ denotes the set of edges leaving set $S$.
The optimization variables are $\routingz \in \{0,1\}^V$, which indicates whether the demand at each node is satisfied, $\vars^k \in \{0,1\}^V$, which indicates the nodes served by vehicle $k$, and $\routingw^k \in \{0,1\}^E$, which indicates the edges traversed by vehicle $k$.
The objective is to maximize the total demand served across all nodes.
Constraint~\eqref{eq:visitnodeonce} ensures each node is visited by at most one vehicle, while constraint~\eqref{eq:capacity} ensures the total demand served by each vehicle does not exceed its capacity~$\gamma_k$.
Constraint~\eqref{eq:vrp:3} ensures that if vehicle $k$ serves node $v$, it must visit that node.
Constraint~\eqref{eq:vrp:start} ensures each vehicle starts at the depot, and constraint~\eqref{eq:vrp:flow} enforces flow conservation so that each vehicle entering a node must also leave it.
The subtour elimination constraints~\eqref{eq:vrp:subtour} prevent vehicles from forming cycles that exclude the depot.
Constraint~\eqref{eq:vrp:fuel} ensures each vehicle does not exceed its fuel limit, and constraint~\eqref{eq:vrp:reward} sets $\routingz_v = 1$ if any vehicle serves node $v$.

We decompose this problem into upper and lower-level subproblems.
In the upper-level problem, we assign nodes to vehicles by solving
\begin{equation*} \ifarxiv \tag{$\mathcal{\hat{P}}_{\text{VRP}}^{1}$} \fi
    \begin{array}{ll}
        \text{maximize} & \predcost(\param)^T \vars \\
        \text{subject to} & \sum_{k=1}^m \vars^k_{v} \le 1, \quad  v \in V \setminus \{v_0\},  \\ & \sum_{v \in V} \demand_v(\param) \vars^k_{v} \le \gamma_k, \quad k=1, \dots, m, \\ 
                          & \vars^{k} \in \{0, 1\}^V, \quad \vars = (\vars^1, \dots, \vars^m),
    \end{array}
\end{equation*}
where $\predcost_\nnparam(\param) \in \reals^{|V| \times m}$ is a learned cost vector and $\vars^k_{v} = 1$ if vehicle $k$ serves node $v$.
The optimization variable is $\vars = (\vars_1, \dots, \vars_m) \in \{0,1\}^{|V| \times m}$.
Once each $\vars^k$ is fixed, the $m$ lower-level problems correspond to single-vehicle routing problems for each $k=1, \dots, m$ of the form,
\begin{equation*}  \ifarxiv \tag{$\mathcal{P}_{\text{VRP}}^{2}$} \fi
    \begin{array}{ll}
        \text{maximize} & \sum_{v \in V} \demand_v(\param) \routingz_v \\
        \text{subject to}
        & \sum_{v \in V} \demand_v(\param) \vars^k_v \le \gamma_k,\\
        & \vars^{k}_{v} \le \sum_{e \in E_v} \routingw^{k}_{e}, \quad  v \in V,\\
        & A \routingw^{k} = 0,\\
        & \vars^{k}_{v} \le \sum_{u \in E_S}  \routingw^{k}_{u}, \quad v \in S, v_0 \notin S, \\
        & \sum_{e \in E} l_e \routingw^{k}_{e} \le f_k, \quad  e \in E, \\
        & \routingz_v \le \sum_{k=1}^m \vars^{k}_{v}, \quad  v \in V,\\
        & \routingz \in \{0, 1\}^V, \quad \routingw^k \in \{0, 1\}^E.
    \end{array}
\end{equation*}
The optimization variables in the lower-level problem are $y = (\routingw_1, \dots, \routingw_m, \routingz) \in \integers^{m|E| + |V|}$.
Note that there are exponentially many constraints in~\eqref{eq:vrp}, as is standard in VRP formulations.
During evaluation we add such constraints lazily with callbacks.

\paragraph{Problem generation.}
We assume all problem parameters are fixed except for the demand at each node $v$, denoted $\delta_v(\param)$, which depends on $\param$.
We consider a complete graph with $\ROUTINGNNODE$ nodes and $\ROUTINGNVEHICLE$ heterogeneous vehicles.
The edge lengths $l_e$ are drawn uniformly from $[0,1]$ and held fixed.
Each vehicle's capacity is drawn uniformly from $[0,\ROUTINGCAPACITY]$, and each vehicle's fuel is drawn uniformly from $[0,\ROUTINGFUEL]$.
To construct the parametric family, we generate a matrix $A \in \reals^{\ROUTINGNNODE \times p}$ from a Gaussian distribution and define the demand vector as $\delta(\param) = |A \param| \in \reals^{\ROUTINGNNODE}$.
The parameter values $\param$ are drawn independently from a $p$-dimensional Gaussian distribution.

\begin{figure}
    \centering
    \includegraphics[scale=0.8, trim={0 0.3cm 0 0}]{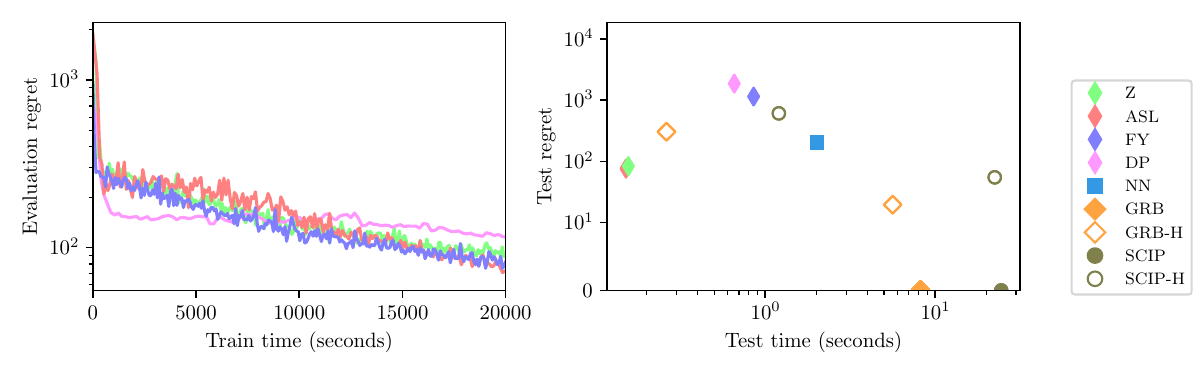}
    \caption{Results for the VRP experiment.
    The left plot shows evaluation regret throughout the training process.
     The right plot shows the test regret~\eqref{eq:absgap} and time (averaged over the test instances) on the test set.
     The first heuristic to terimnate for Gurobi and SCIP terminates after finding a single feasible solution. The second to terminate terminates after finding three.}
    \label{fig:vrp}
\end{figure}

\paragraph{Results.}
Results for the experiment are plotted in Figure~\ref{fig:vrp}.
Each of our trained models perform better than nearest neighbor, and they all find solutions significantly faster than Gurobi.
All models find feasible solutions faster than Gurobi or SCIP do, and the feasible solutions found are of higher quality than the first solutions found by Gurobi.

\begin{table}[ht]
\centering
{\fontsize{7}{8.5}\selectfont
\begin{tabular}{l*{9}{r}}
\toprule
 & NN & DP & GRB & GRB-H & SCIP & SCIP-H & \textbf{ASL} & \textbf{Z} & \textbf{FY} \\
\midrule
  \begin{filecontents*}{routing_scientific.csv}
name,NN,DP,GRB,GRB_H,SCIP,SCIP_H,ASL,Z,FY
Time (s),2.02e+00,6.60e-01,8.22e+00,5.63e+00,2.45e+01,2.95e+01,1.50e-01,1.60e-01,8.50e-01
Regret,2e+02,1.9e+03,1.00e-01,1.96e+01,0.00e+00,5.48e+01,7.62e+01,8.35e+01,1.2e+03
Normalized Regret,1.89e+00,1.65e+01,8.00e-04,1.74e-01,0.00e+00,5.28e-01,6.84e-01,7.41e-01,1.01e+01
\end{filecontents*}
\csvreader[
      head to column names,
      late after line=\\,
      before reading={\catcode`\ =9\relax}
      ]{routing_scientific.csv}{}{%
      \csvcoli & \csvcolii & \csvcoliii & \csvcoliv & \csvcolv &
      \csvcolvi & \csvcolvii & \csvcolviii & \csvcolix & \csvcolx
      }
    \bottomrule
\end{tabular}
}
\caption{Results for the routing experiment.
Our method names are in bold.
The heuristic methods recorded here terminate after finding a single feasible solution.}
\end{table}

\subsection{Conformal prediction}
\begin{figure}
    \centering
    \includegraphics[scale=0.75]{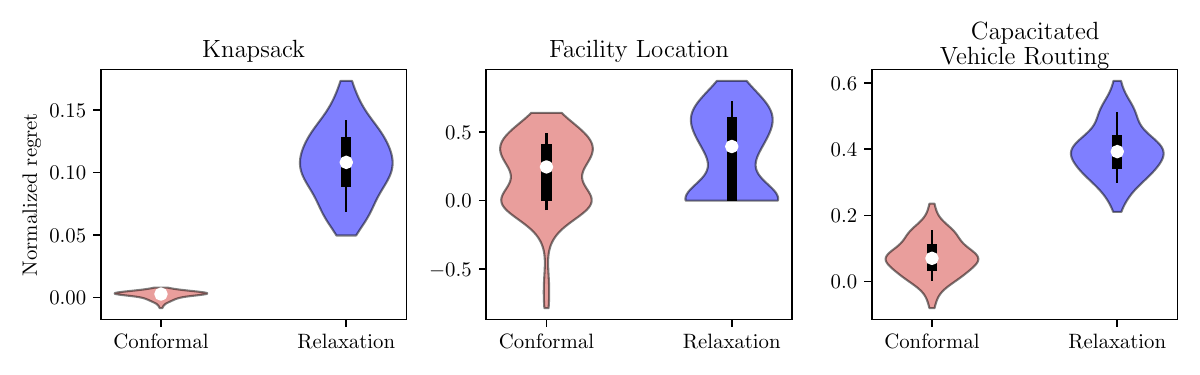}
    \caption{Results for the conformal prediction experiments.
    In each plot the bound given by the learned conformal prediction model is compared to the relaxation bound.
The density plot is over the test set which is unseen during training.} \label{fig:conformal}
\end{figure}

We fit the conformal predictor to the best-performing model in each case, selected based on the lowest evaluation regret.
The conformal prediction model is fitted as described in Section~\ref{sec:conformal}: The bound prediction function $\confnn_\confparam$ is trained on the evaluation set $\eval$ and calibrated on the calibration set $\calibration$.
We evaluate the quality of the bounds on the test dataset $\test$.
In Figure~\ref{fig:conformal}, we compare the normalized regret $(\xi(\param) - z(\param))/|z(\param)|$ for two bound functions $\xi(\param)$.
The first, labeled \textit{conformal}, outputs the bound given by the conformal prediction model.
The second, labeled \emph{relaxation}, is the convex relaxation bound for the problem parametrized by $\param$.
We display density plots of the bound values for each problem in the test dataset.
As shown in Figure~\ref{fig:conformal}, the conformal prediction method yields significantly tighter bounds than the relaxation.
However, the conformal bounds may be invalid up to approximately an $\alpha$ fraction of the time, which occurs when the normalized bound becomes negative.
Results for the tightness of the predicted conformal bounds are summarized in Table~\ref{tab:conformal}.
\begin{table}[ht]
\centering
{\fontsize{7}{8.5}\selectfont
\begin{tabular}{l*{8}{r}}
\toprule
 & \multicolumn{2}{c}{Knapsack} & \multicolumn{2}{c}{FacLoc} & \multicolumn{2}{c}{VRP} \\
 Method & Conformal & Relaxation & Conformal & Relaxation & Conformal & Relaxation \\
 \midrule
\begin{filecontents*}{full_conformal_results_tables_paper_1.csv}
Problem,Knapsack,Knapsack,Facility Location,Facility Location,Capacitated Vehicle Routing,Capacitated Vehicle Routing,Stochastic Vehicle Scheduling,Stochastic Vehicle Scheduling
Regret,2.9e+03,1.4e+05,4e+03,7.2e+03,8.72,46.25,0.33,5.03
Percent wrong,0.21,0.00,0.24,0.00,0.08,0.00,0.10,0.00
Normalized Regret (correct),0.03,0.00,0.54,0.00,0.20,0.00,0.07,0.00
Normalized Regret (incorrect),0.03,0.00,1.9e+03,0.00,0.09,0.00,0.03,0.00
\end{filecontents*}
\csvreader[
      head to column names,
      late after line=\\,
      before reading={\catcode`\ =9\relax}
      ]{full_conformal_results_tables_paper_1.csv}{}{%
      \csvcoli & \csvcolii & \csvcoliii & \csvcoliv & \csvcolv &
      \csvcolvi & \csvcolvii}
\bottomrule
\end{tabular}
}
\caption{Results for the bounds given by conformal prediction.} \label{tab:conformal}
\end{table}

\section{Conclusion}
We presented a hierarchical learning architecture for efficiently computing high-quality solutions to structured mixed-integer programs.
Our approach decomposes the original problem into smaller upper and lower-level subproblems that are solved sequentially, with the upper-level decisions parametrizing the lower-level constraints.
We formulated the training problem as a convex optimization task using decision-focused learning techniques, developing several surrogate convex losses that approximate the true loss while providing meaningful gradients and preserving optimality.
To provide robustness guarantees, we introduced a conformal prediction method that yields probabilistic bounds on solution suboptimality by training a predictor on an evaluation dataset and calibrating it on a separate calibration dataset.
Numerical experiments on facility location, knapsack, and vehicle routing problems demonstrated that our approach finds high-quality feasible solutions significantly faster than state-of-the-art MIP solvers.
While the problems we tested are relatively small (requiring 10 to 100 seconds to solve with standard solvers), our results demonstrate the viability of learned approaches for hierarchical MIPs.
The conformal prediction techniques we developed are particularly promising, as they can be applied to other feasible-solution-finding heuristics for MIPs.
Future work should focus on scaling these methods to larger problem instances and exploring their application to broader classes of structured optimization problems.

\ifarxiv
\bibliography{bibliography}
\fi

\end{document}